\documentclass[12pt]{amsart}
\usepackage{amssymb}
\usepackage{epsfig}
\newcommand{\Lambdabar}{\overline{\Lambda}}

\newcommand{\gammabar}{\overline{\gamma}}
\newcommand{\Lbar}{\overline{L}}
\newcommand{\Fbar}{\overline{F}}
\newcommand{\F}{\mathbb{F}}
\newcommand{\R}{\mathbb{R}}
\newcommand{\C}{\mathbb{C}}
\newcommand{\Z}{\mathbb{Z}}
\newcommand{\Q}{\mathbb{Q}}
\newcommand{\Nlt}{N_{\rm left}}
\newcommand{\Nrt}{N_{\rm right}}
\newcommand{\Res}{{\rm Res}}
\newcommand{\Ind}{{\rm Ind}}
\newcommand{\Gal}{{\rm Gal}}
\newcommand{\GL}{{\rm GL}}
\newcommand{\SL}{{\rm SL}}
\newcommand{\PGL}{{\rm PGL}}
\newcommand{\Tr}{{\rm Tr}\,}
\newcommand{\Ad}{{\rm Ad}}
\newcommand{\ord}{{\rm ord}}
\newtheorem{theorem}{Theorem}
\newtheorem{prop}[theorem]{Proposition}
\newtheorem{lemma}[theorem]{Lemma}

\theoremstyle{definition}
\newtheorem{definition}[theorem]{Definition}
\theoremstyle{remark}
\newtheorem*{remark}{Remark}
\numberwithin{equation}{section}
\numberwithin{theorem}{section}
\numberwithin{figure}{section}
\numberwithin{table}{section}
\begin{document}
\title{Artin's conjecture, Turing's method and the Riemann hypothesis}
\date{}
\author{Andrew R.\ Booker}
\thanks{The author was supported by an NSF postdoctoral fellowship}
\address{
Mathematics Department\\
530 Church Street\\
University of Michigan\\
Ann Arbor, MI 48109}
\email{{\tt arbooker@umich.edu}}
\begin{abstract}
We present a group-theoretic criterion under which one may verify the
Artin conjecture for some (non-monomial) Galois representations, up to
finite height in the complex plane.  In particular, the criterion applies
to $S_5$ and $A_5$ representations.  Under more general conditions, the
technique allows for the possibility of verifying the Riemann hypothesis
for Dedekind zeta functions of non-abelian extensions of $\Q$.

In addition, we discuss two methods for locating zeros of arbitrary
$L$-functions.  The first uses the explicit formula and techniques
developed in \cite{bs} for computing with trace formulae.
The second method generalizes that of Turing for
verifying the Riemann hypothesis.  In order to apply it we develop a
rigorous algorithm for computing general $L$-functions on the critical
line via the Fast Fourier Transform.

Finally, we present some numerical results testing Artin's conjecture
for $S_5$ representations, and the Riemann hypothesis for Dedekind zeta
functions of $S_5$ and $A_5$ fields.
\end{abstract}
\maketitle
\section{Introduction}
\subsection{Artin's conjecture}
Let $K/\Q$ be a Galois extension and $\rho:\Gal(K/\Q)\to\GL_n(\C)$ a
non-trivial, irreducible representation of its Galois group.  In
\cite{artin}, Artin associated to this data an $L$-function
$L(s,\rho)$, defined initially for $\Re(s)>1$, which he conjectured
to continue to an entire function and satisfy a functional equation.
By a theorem of Brauer \cite{brauer}, one now knows the {\em
meromorphic} continuation and functional equation of Artin's
$L$-functions.  The question remains whether they can have poles in
the critical strip $0<\Re(s)<1$.

Artin established his conjecture for the {\em monomial}
representations, those induced from a $1$-dimensional representation of
a subgroup; this of course includes all $1$-dimensional $\rho$,
in which case $L(s,\rho)=L(s,\chi)$ for a Dirichlet character
$\chi$.  Although the conjecture has not been decided in any dimension
$\ge 2$, more evidence is provided in dimension $2$
by the Langlands-Tunnell theorem \cite{langlands,tunnell}, which
affirms the conjecture for those representations whose image in
$\PGL_2(\C)$ is isomorphic to $A_4$ (tetrahedral) or $S_4$
(octahedral); only the $A_5$ (icosahedral) case remains.  When
$\rho$ is an {\em odd} icosahedral representation, meaning $\det\rho$
determines an odd Dirichlet character, infinitely many examples of
Artin's conjecture are known from the work of Taylor et
al.\ \cite{taylor1,taylor2}.

Moreover, in the odd $2$-dimensional case, there is an algorithm for
verifying the conjecture, as follows.  By a construction of Deligne
and Serre \cite{deligne-serre}, given a holomorphic modular form $f$ of
weight $1$, one may associate an odd 2-dimensional representation $\rho$
such that $L(s,f)=L(s,\rho)$.  Conversely, every odd $2$-dimensional
$\rho$ such that $L(s,\rho)$ is entire arises from the Deligne-Serre
construction.  For any particular $\rho$, one can search for the
associated form; once found, comparing the representation constructed
by Deligne-Serre to $\rho$ via an effective version of the Cebotarev
density theorem allows one to deduce the conjecture for $\rho$.  This and
other related techniques have been carried out in a number of cases;
see \cite{buhler,kiming,jehanne,stein}.

On the other hand, if one considers {\em even} $2$-dimensional
representations, the situation is somewhat different.  There as well
the conjecture has been established for all but the icosahedral cases.
However, the correspondence is not with holomorphic forms, but rather
Maass forms of eigenvalue $\frac14$.  Unfortunately, no analogue of the
result of Deligne and Serre is known in that setting.  Moreover,
computation of the associated forms remains elusive; existing
techniques (see e.g.\ \cite{bsv}) only allow one to calculate Maass
forms to within a prescribed precision, never exactly.  Thus, at
present this approach does not yield an algorithm for verifying Artin's
conjecture.

The apparent difference between these two cases leads naturally to the
following question:  Given a Galois representation $\rho$, is there an
algorithm that will decide in finite time, with proof, whether
$L(s,\rho)$ is entire?  Note that like the Riemann hypothesis, Artin's
conjecture is falsifiable, i.e.\ it may be disproven by observing a
counterexample, in this case a pole.  The challenge is thus to find a
way of demonstrating the conjecture when true.

Although we are unable to provide a definitive answer to this question,
one approach, at least for $2$-dimensional
representations, is suggested by a theorem from \cite{booker}: If a
given $2$-dimensional $\rho$ is not associated to a holomorphic or
Maass form as above, then $L(s,\rho)$ has infinitely many poles.  In
particular, once $L(s,\rho)$ has at least one pole, it must have
infinitely many.  Unfortunately, the result is ineffective, in the
sense that it does not predict where the first pole must occur.  A
natural question, therefore, is whether an effective version of this
theorem exists.  First, however, we must consider exactly what that
would mean; since the only handle that we have on an Artin $L$-function
in the critical strip is as the ratio of entire functions given by
Brauer's theorem, it is not immediately clear that we can check its
holomorphy at a zero of the denominator without {\em a priori} knowing a
lower bound on the residue of any poles.

In this paper we address precisely this issue, in Section \ref{sec:artin}.
There we present a criterion which, when satisfied, yields an algorithm
for verifying the holomorphy of an Artin $L$-function up to a given
height in the critical strip.  In particular, we give the first direct
evidence (as far as we are aware) of holomorphy in the critical strip
of an $L$-function for which the conjecture cannot be established
through the methods mentioned above.  Although our criterion is not
always satisfied, we are in general able to deduce partial information,
such as a bound on the multiplicities and residues of possible poles.
Moreover, the limitations of the information that we obtain give an idea
of the hypotheses that one would have to impose in any effective version
of the converse theorem in order to make the above approach work.

\subsection{Turing's method and the Riemann hypothesis}
One application of our criterion is to the Riemann hypothesis for
Dedekind zeta functions.  Turing \cite{turing} devised a method for
checking the hypothesis in a bounded region for the Riemann $\zeta$
function.\footnote{Reading Turing's paper on the subject, which was his
last, one marvels at what he accomplished with the limited computational
resources of the day.  His method was truly ahead of its time.}  The
method depends on the simplicity of the zeros of $\zeta$.  Because of
that, it is only directly extendable to Dedekind zeta functions of
non-normal extensions of small degree (see \cite{tollis}) or abelian
extensions, for which it is more natural to verify the hypothesis for
the associated Dirichlet $L$-functions instead (see \cite{rumely}).

Similarly, for a non-abelian extension one can factorize the zeta
function into Artin $L$-functions of irreducible representations.
As these are also expected to have simple zeros, Turing's method
applies, provided one assumes the Artin conjecture.  However, combining
our criterion with Turing's method, we will in some cases be able to
deduce the Riemann hypothesis and holomorphy of the relevant Artin
$L$-functions simultaneously.  In fact, as we will see, there are even
cases where we may check the Riemann hypothesis without being able to
verify Artin's conjecture.  We carry out the necessary generalization
of Turing's method in Section \ref{turing}.

\subsection{Rigorous zero computations}
In order to implement these ideas, we develop, in Sections
\ref{sec:explicit} and \ref{sec:rigorous}, two methods of locating zeros
of $L$-functions.  The first uses the explicit formula and techniques
developed for the Selberg trace formula in \cite{bs}.  If one assumes the
Riemann hypothesis, this method may be used with our criterion, in place
of Turing's method, for verifying the Artin conjecture.  More importantly,
the explicit formula is clean to implement and yields estimates for low
zeros quickly.  It can thus serve as a check for later computations,
or to fine tune the parameters of Turing's method for greater speed.

The second method is a technique for fast, rigorous computations of
$L$-functions on the critical line.  This is a hard problem in general,
basically because of the difficulty of providing uniform, effective
bounds for the relevant Mellin transforms.  By making use of the Fast
Fourier Transform, our technique allows one to compute many values of
the same $L$-function simultaneously, which is particularly appropriate
for Turing's method.  In doing so, we need only consider a single Mellin
transform, making rigorous computation more practical.  In addition,
the method has complexity comparable to that of computing a {\em single}
value by the approximate functional equation.

Although our primary interest is in Artin $L$-functions, we carry out
the details of Sections \ref{sec:explicit}, \ref{turing} and
\ref{sec:rigorous} for {\em arbitrary} $L$-functions $L(s)$, in the
hope that the results may be useful outside of the present context.
More precisely, we make the following assumptions, notations and
conventions throughout:
\begin{itemize}
\item $L(s)$ is given by an Euler product of degree $r$:
\begin{equation}
L(s)=\prod_{p\,\mbox{\scriptsize prime}}
\frac1{(1-\alpha_{p,1}p^{-s})\cdots(1-\alpha_{p,r}p^{-s})},
\label{eulerprod}
\end{equation}
where the $\alpha_{p,j}$ are complex parameters satisfying
the individual bound $|\alpha_{p,j}|\le p^{1/2}$, and the product is
absolutely convergent for $\Re(s)>1$.  Further, for all but finitely
many $p$, there is a pairing $\alpha\mapsto\alpha'$ such that
$|\alpha_{p,j}\alpha_{p,j}'|=1$.  For the exceptional $p$, such a
pairing exists for a subset of the $\alpha_{p,j}$, and those not in the
subset satisfy $|\alpha_{p,j}|\le 1$.
\item Define
\begin{equation}
\Gamma_{\R}(s) := \pi^{-s/2}\Gamma\!\left(\frac{s}2\right),
\qquad
\gamma(s) := \epsilon N^{\frac12(s-\frac12)}\prod_{j=1}^r\Gamma_{\R}(s+\mu_j),
\qquad \Lambda(s) := \gamma(s)L(s),
\end{equation}
where $|\epsilon|=1$, $N$ is a positive integer and $\Re(\mu_j)\ge -\frac12$.
For a certain choice of these parameters, $\Lambda(s)$ has meromorphic
continuation to $\C$, is a ratio of entire functions of order $1$,
and satisfies the functional equation
\begin{equation}
\Lambda(s) = \Lambdabar(1-s),
\label{funceq}
\end{equation}
where for a complex function $f$ we denote by $\overline{f}(s)$ the
function $\overline{f(\bar{s})}$.  Note that $\epsilon$ here is the
square root of the usual root number, and is only defined up to
multiplication by $\pm1$; we choose the value with argument in
$\bigl(-\frac{\pi}2,\frac{\pi}2\bigr]$.  Including $\epsilon$ as part
of the $\gamma$ factor makes $\Lambda(s)$ real for $\Re(s)=\frac12$, as
can be seen from \eqref{funceq}.

Let $Q(s)$ be the {\em analytic conductor}:
\begin{equation}
Q(s):=N\prod_{j=1}^r\frac{s+\mu_j}{2\pi}.
\end{equation}
Note that $\gamma(s)$ satisfies the recurrence
$\gamma(s+2)=Q(s)\gamma(s)$.  Further, we define
\begin{equation}
\chi(s):=\frac{\gammabar(1-s)}{\gamma(s)},
\end{equation}
so that $L(s)=\chi(s)\Lbar(1-s)$.
\item $L(s)$ may have at most finitely many poles, which we assume
to lie along the line $\Re(s)=1$.  We label them $1+\lambda_k$ with
$\lambda_k\in i\R$, $k=1,\ldots,m$, repeating with the appropriate
multiplicity.  Further, from the functional
equation \eqref{funceq}, each $\lambda_k$ will equal $-\mu_j$ for some
$j$, counting multiplicity; in particular, $m\le r$.  We set
\begin{equation}
P(s) := \prod_{k=1}^m(s-\lambda_k),
\end{equation}
so that $P(s)P(s-1)\Lambda(s)$ is entire.
\item Some progress is known toward the
Ramanujan conjecture for $L$; that is, there exists $\theta<\frac12$
such that
\begin{equation}
|\alpha_{p,j}| \le p^{\theta}
\qquad\mbox{and}\qquad
\Re(\mu_j)\ge -\theta
\label{ramanujan}
\end{equation}
for all $p,j$.  This assumption is not strictly necessary, as we could
instead use average bounds of Rankin-Selberg type.  However, bounds of
the form \eqref{ramanujan} are now known in the cases of greatest
interest (automorphic $L$-functions \cite{lrs}), and the results are
easier to state and use assuming it.
\end{itemize}

\subsection{Numerical results}
Finally, in Section \ref{numerics} we describe the implementation
of the above ideas and give some numerical results of tests of the
Riemann hypothesis for a few $S_5$ and $A_5$ extensions in the region
$|\Im(s)|\le 100$.  For the $S_5$ cases, this includes a verification
of Artin's conjecture in the same region for the $L$-functions of all
representations of the group.

\section*{Acknowledgements}
The idea for this paper arose from a conversation with Laurent Clozel,
and it is a pleasure to acknowledge him here.  I thank Harold Stark, who
has considered problems of this nature in the past and whose interest
motivated me to finish the project.  Thanks also to Bob Griess, Martin
Isaacs, Peter Sarnak and Kannan Soundararajan for helpful discussions.

\section{A criterion for verifying Artin's conjecture}
\label{sec:artin}
Let $\rho:\Gal(K/\Q)\to\GL_n(\C)$ be a Galois representation, as in the
introduction.  Brauer's theorem expresses the $L$-function $L(s,\rho)$
as a ratio $N(s)/D(s)$, where $N(s)$ and $D(s)$ are Artin $L$-functions
associated to sums of monomial representations.  If $f(s)$ is any
holomorphic Artin $L$-function, we have a formula for the number
$N_f(t_1,t_2)$ of zeros of $f$ between heights $t_1$ and $t_2$, from the
argument principle:
\begin{equation}
N_f(t_1,t_2)=\frac1{2\pi i}\int_C\frac{f'}{f}(s)\,ds,
\label{eq:argp1}
\end{equation}
where $C$ is the rectangular contour with vertices at $2+it_1$,
$2+it_2$, $-1+it_2$, $-1+it_1$ and counter-clockwise orientation.
We also have available in this case algorithms to compute $f$ and $f'$ at
an arbitrary point in the complex plane; see Section \ref{sec:rigorous}.
Thus, in principle we could compute \eqref{eq:argp1} exactly by numerical
integration.  Although \eqref{eq:argp1} has the advantage of applying in
great generality, to do so would be inefficient and difficult to implement
rigorously.  In the special case that the zeros of $f$ are simple, a much
more efficient algorithm was given by Turing; see Section \ref{turing}.

No matter how we arrive at the numbers $N_f(t_1,t_2)$, there is always
some uncertainty in the locations of the zeros of $f$.  In
\eqref{eq:argp1} this is due to the fact that as $t_i$ approaches the
ordinate of a zero, higher and higher precision is needed in order to
compute $f'/f$ accurately.  This is in line with the expectation that
the general zero is transcendental, meaning that one can never know it
exactly.

For $L(s,\rho)$, we can recover the {\em net} number of zeros
(i.e.\ zeros minus poles) between heights $t_1$ and $t_2$ as
$N_N(t_1,t_2)-N_D(t_1,t_2)$.  If Artin's conjecture is true then for
every zero of $D(s)$ there is a zero of $N(s)$ at the same point.
However, because of the uncertainty in the locations of the zeros of
$N(s)$ and $D(s)$, from this computation alone we cannot rule out the
possibility that $L(s,\rho)$ has a pole with a zero very close by in
the neighborhood of a zero of $D(s)$.  In other words, we can only
observe the counts of net zeros in these small neighborhoods.

Fortunately, there are some restrictions on potential poles.  For
instance, the Dedekind zeta function of the extension, $\zeta_K(s)$,
factors into Artin $L$-functions:
\begin{equation}
\zeta_K(s) = \prod_{\rho}L(s,\rho)^{\dim\rho},
\end{equation}
where the product is over all irreducible representations of
$\Gal(K/\Q)$.  Since $\zeta_K(s)$ is holomorphic (except for a simple
pole at $s=1$), we see that any pole of $L(s,\rho)$ in the critical
strip must be located at the zero of another function.
More generally, if $\sigma$ is any representation, we have
\begin{equation}
L(s,\sigma) = \prod_{\rho}L(s,\rho)^{\langle\sigma,\rho\rangle},
\end{equation}
where $\langle\cdot,\cdot\rangle$ is the inner product on the space of
characters, and by abuse of notation we write
$\langle\sigma,\rho\rangle$ for $\langle\Tr\sigma,\Tr\rho\rangle$.
When $\sigma$ is monomial, we again have $L(s,\sigma)$ holomorphic with
the possible exception of a pole at $s=1$.

This information is described most concisely by use of the {\em
Heilbronn (virtual) character}:  For $s_0\in\C\setminus\{1\}$, define
\begin{equation}
\theta_{s_0} = \sum_{\rho}\ord_{s=s_0}L(s,\rho)\cdot\Tr\rho,
\end{equation}
where $\ord_{s=s_0}L(s,\rho):=\Res_{s=s_0}\frac{L'}{L}(s,\rho)$.
Thus,
\begin{equation}
\ord_{s=s_0}L(s,\sigma)=\langle\theta_{s_0},\sigma\rangle
\ge 0\quad\mbox{for all monomial }\sigma.
\label{eq:monoholo}
\end{equation}
The study of Heilbronn characters leads to many useful results.  For
example, in \cite{foote-murty} it is shown that
\begin{equation}
\sum_{\rho}\bigl(\ord_{s=s_0}L(s,\rho)\bigr)^2
\le \bigl(\ord_{s=s_0}\zeta_K(s)\bigr)^2.
\label{eq:zetaKbound}
\end{equation}
In particular, the zeros and poles of each $L(s,\rho)$ are
among the zeros of $\zeta_K(s)$.

The idea now is to combine \eqref{eq:monoholo} with observations of
net zeros.  If we look in a small enough neighborhood of a zero of
$\zeta_K(s)$, we expect to find one net zero for a single $L(s,\rho)$
and no net zeros for the others.  This is based on the assumption that
the zeros of different irreducible Artin $L$-functions are distinct and
simple.  While such a statement is likely impossible to prove, we may use
it as a working hypothesis to be tested at run time.  This is analogous
to assuming the simplicity of the zeros of $\zeta$ in order to check the
Riemann hypothesis.  (Note that if there were a multiple zero of $\zeta$,
it is doubtful that one could distinguish it from a counterexample.)

In other words, if the working hypothesis is true, then our net zero
observations correspond to the character $\Tr\rho$ for some $\rho$.
Thus, we have $\Tr\rho=\theta_{s_1}+\ldots+\theta_{s_n}$, where
$s_1,\ldots,s_n$ are the distinct zeros of $\zeta_K(s)$ in the
neighborhood that we examine.  We would like to conclude that there is
just one such point, meaning that the actual zero counts agree with our
observations.  Since the Heilbronn characters satisfy
\eqref{eq:monoholo}, it is enough to show that
\begin{equation}
\Tr\rho\ne\chi_1+\chi_2\mbox{ for virtual characters }
\chi_i\ne 0\mbox{ with }\langle\chi_i,\sigma\rangle\ge 0
\mbox{ for all monomial }\sigma.
\label{eq:monohyp}
\end{equation}

The one notable exception to this philosophy is at the central point
$\frac12$, where there can be forced vanishing if $\rho$ is self-dual
(an example of which is given in \cite{armitage}).
In that case, we expect one zero for each self-dual $\rho$ with an odd
functional equation, and no zeros for the rest.  However, we can only
determine the parity of the order of vanishing at $\frac12$.  This leads
to the following replacement for condition \eqref{eq:monohyp} at $\frac12$:
\begin{equation}
\sum_{\substack{\rho\;{\rm self-dual}\\ \Lambda(1-s,\rho)=-\Lambda(s,\rho)}}
\Tr\rho\ne\chi_1+2\chi_2\quad\mbox{with }\chi_i\ne 0\mbox{ and }
\langle\chi_i,\sigma\rangle\ge 0\mbox{ for all monomial }\sigma.
\label{eq:monohyp2}
\end{equation}

\begin{figure}
\centerline{\epsfxsize=3in \epsfbox{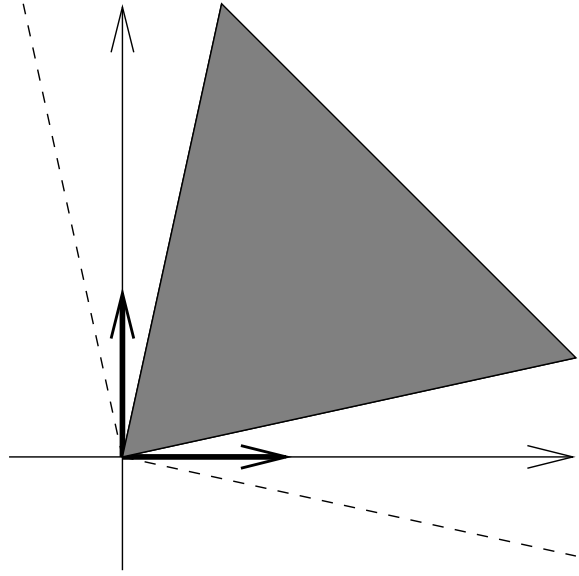}}
\caption{Lattice of virtual characters.  The shaded cone is the monoid
generated by monomial characters.}
\label{lattice}
\end{figure}

When \eqref{eq:monohyp} is satisfied for all irreducible representations
$\rho$, we may check the holomorphy of all $L(s,\rho)$ at any point
at which the working hypothesis turns out to be true.  We give a name to
describe this situation:
\begin{definition} A finite group $G$ is {\em almost monomial }if, for
each irreducible representation $\rho$, if $\Tr\rho = \chi_1+\chi_2$
for virtual characters $\chi_i$ such that
$\langle\chi_i,\sigma\rangle\ge 0$ for all monomial $\sigma$, then
either $\chi_1 = 0$ or $\chi_2$ = 0.
\end{definition}
The terminology is explained with aid of Figure \ref{lattice}.  The
plane represents the lattice of virtual characters, with the first
quadrant being the monoid of characters, and the shaded cone
the monoid generated by the monomial characters.  We consider all
virtual characters within $90$ degrees of the cone, which in the
figure is everything within the dashed lines.  The group is almost
monomial if this set is not much larger than the character monoid, in
the precise sense that the irreducible representations, which are the
coordinate axis vectors represented by thick arrows, remain
indecomposable in this set.  Equivalently, the monoid generated by the
monomial characters should be close to the full character monoid.
From the picture it is easy to see that any monomial group is almost
monomial.

One could argue that we should include condition \eqref{eq:monohyp2} in
our definition as well.  We prefer to keep it separate, taking the view
that it is more important to be able to demonstrate holomorphy at a
generic zero of the denominator.  Indeed, we have already seen that the
$L$-function of a $2$-dimensional representation cannot have a finite
number of poles, so we do not lose much generality by excluding a
single point.  It is plausible that such a result holds for higher
dimensions as well.  Moreover, condition \eqref{eq:monohyp2} seems
usually to be weaker than almost monomiality; cf.\ Proposition
\ref{prop:monoex} below.

A potentially more serious issue is that $N(s)$ and $D(s)$ may have
high order zeros at $\frac12$, in which case Turing's method does not
apply.  This could be remedied by computing the contour integral
\eqref{eq:argp1} around $\frac12$, but we would like to avoid doing
so.  Fortunately, if the order of vanishing at $\frac12$ is at most
$3$, we can still conclude that we have the correct count by sign
changes alone; that is because for a self-dual representation, if we
miss a zero away from $\frac12$ then we must miss at least four such
zeros.  Fortunately again, in all cases that we consider, $N(s)$ has at
most three irreducible factors with a potential zero at $\frac12$.

Like monomiality, the notion of almost monomiality behaves well under
some common group operations.  In particular, we have the following.
\begin{prop}
If $G$ is almost monomial then so are quotients of $G$ and products
$G\times H$ for any monomial group $H$.
\label{prop:monofunc}
\end{prop}
\begin{proof}
1. Let $K$ be a normal subgroup of $G$ and $\tilde{\pi}$ an irreducible
representation of $G/K$.  Suppose that
$\Tr\tilde{\pi}=\tilde{\chi}_1+\tilde{\chi}_2$, with
$\langle\tilde{\chi}_i,\tilde{\sigma}\rangle\ge 0$ for all monomial
$\tilde{\sigma}$.  Let $\pi$, $\chi_i$ be the lifts of
$\tilde{\pi}$, $\tilde{\chi}_i$ to $G$ obtained by composition with the
natural projection.  Then $\pi$ is irreducible and $\Tr\pi=\chi_1+\chi_2$.
Further, if $\rho$ is an irreducible representation of $G$ then
$\langle\chi_i,\rho\rangle=0$ unless $\rho$ factors through $G/K$.  If
that is the case, let $\tilde\rho$ denote the induced map on $G/K$.

Now, if $\sigma=\Ind_H^G\lambda$ is a monomial representation then, by
Frobenius reciprocity, we have
$\langle\sigma,\rho\rangle=
\bigl\langle\Res_H^G\rho,\lambda\bigr\rangle$
for all $\rho$ factoring through $G/K$.  If at least one of these is
non-zero, i.e.\ $\lambda$ occurs in $\Res_H^G\rho$, then since $\rho$
factors through $G/K$, $\lambda$ must factor through $H/H\cap K
\cong HK/K$.  Let $\tilde{\lambda}$ denote the induced map on $HK/K$.
Then
$\bigl\langle\Res_{HK/K}^{G/K}\tilde{\rho},\tilde{\lambda}\bigr\rangle=
\bigl\langle\Res_H^G\rho,\lambda\bigr\rangle$.
Thus, $\tilde{\sigma}=\Ind_{HK/K}^{G/K}\tilde{\lambda}$ satisfies
$\langle\sigma,\rho\rangle=
\langle\tilde{\sigma},\tilde{\rho}\rangle$.
Therefore,
$\langle\chi_i,\sigma\rangle
=\langle\tilde{\chi}_i,\tilde{\sigma}\rangle\ge 0$.
The conclusion follows by almost monomiality of $G$.

2. Let $\rho_G$ and $\rho_H$ be irreducible representations of $G$ and
$H$, respectively, and suppose that $\Tr\rho_G\otimes\rho_H
=(\Tr\rho_G)(\Tr\rho_H)=\chi_1+\chi_2$ with
$\langle\chi_i,\sigma\rangle_{G\times H}\ge 0$
for all monomial $\sigma$.  Taking
the inner product over $H$ with $\rho_H$, we get
$\Tr\rho_G=\langle\chi_1,\rho_H\rangle_H+\langle\chi_2,\rho_H\rangle_H$.

Next, if $\sigma_G$ is any monomial representation of $G$, we have
$\bigl\langle\langle\chi_i,\rho_H\rangle_H,\sigma_G\bigr\rangle_G
=\langle\chi_i,\sigma_G\otimes\rho_H\rangle_{G\times H}
\ge 0$, since $\sigma_G\otimes\rho_H$ is monomial.  Thus, since $G$ is
almost monomial, we have $\langle\chi_i,\rho_H\rangle=0$ for some $i$.

Similarly, if $\rho_H'$ is any other irreducible representation of $H$,
we find
$0=\langle\chi_1,\rho_H'\rangle_H+\langle\chi_2,\rho_H'\rangle_H$.
Thus, $\langle\chi_1,\rho_H'\rangle=\langle\chi_2,\rho_H'\rangle=0$.
Therefore $\chi_i=0$ for some $i$.
\end{proof}

The next proposition shows that the class of almost monomial groups is
strictly larger than that of monomial groups.
\begin{prop}
The groups $\SL_2(\F_3)$, $A_5$ and $S_5$ are almost monomial and satisfy
\eqref{eq:monohyp2}.
\label{prop:monoex}
\end{prop}
\begin{proof}
These are shown with the aid of the computer algebra system GAP
\cite{GAP}.  We illustrate the general procedure for checking almost
monomiality for a given group with the example $A_5$.  Note first that
$A_5$ has five irreducible representations, of dimensions $1$, $3$,
$3$, $4$ and $5$.  We use GAP to determine all monomial representations.
In this case they are non-negative linear combinations of the vectors
$(1,0,0,0,0)$,
$(0,0,0,0,1)$,
$(1,0,0,1,0)$,
$(0,1,1,0,0)$,
$(0,1,0,1,1)$,
$(0,0,1,1,1)$,
and
$(0,1,1,1,0)$,
where the components indicate the multiplicities of the irreducible
representations.  We label the monomial representations associated to
these vectors $\sigma_1,\ldots,\sigma_7$.  The first five form
a $\Z$-basis for the lattice of virtual characters, i.e.\ any virtual
character $\chi$ may be written uniquely as an integral linear
combination $\chi=\sum_{i=1}^5x_i\Tr\sigma_i$.

Now, almost monomiality is equivalent to the assertion that for each
irreducible representation $\rho$, whenever
$0\le\langle\chi,\sigma\rangle\le\langle\rho,\sigma\rangle$ for all
monomial $\sigma$, we have either $\chi=0$ or $\chi=\Tr\rho$.  Using
our integral basis, we investigate the solutions to
\begin{equation}
0\le\sum_{i=1}^5x_i\langle\sigma_i,\sigma_j\rangle
\le\langle\rho,\sigma_j\rangle
\label{eq:A5proof}
\end{equation}
for $j=1,\ldots,7$.  Restricting to $j=1,\ldots,5$, we get an invertible
system, i.e.\ the matrix
$A=\bigl(\langle\sigma_i,\sigma_j\rangle\bigr)_{1\le i,j\le5}$ lies in
$\SL_5(\Z)$.  We consider the vectors ${\bf x}=A^{-1}{\bf y}$ for all
${\bf y}=(y_1,\ldots,y_5)$ satisfying
$0\le y_j\le\langle\rho,\sigma_j\rangle$.  By construction, these
satisfy \eqref{eq:A5proof} for $j=1,\ldots,5$.  We check that the only
${\bf x}$ satisfying \eqref{eq:A5proof} for $j=6,7$ are $0$ and
$A^{-1}\bigl(\langle\rho,\sigma_j\rangle\bigr)$, corresponding to $\chi=0$ and
$\chi=\Tr\rho$, respectively.

Similarly, for \eqref{eq:monohyp2} we try all possible combinations
of $\rho$ having odd functional equation.  We may exclude those whose
$L$-functions may be expressed in terms of Dedekind zeta functions, for
which the root number is always $1$.  For $A_5$, the only non-trivial
possibility is that the two $3$-dimensional representations have odd
functional equation.
\end{proof}

With the evidence provided by Propositions \ref{prop:monofunc} and
\ref{prop:monoex}, one might hope that all groups are almost monomial.
That is not the case, as the counterexamples $\GL_2(\F_3)$ and
$\SL_2(\F_5)$ show.  $\SL_2(\F_5)$ has irreducible representations of
dimensions $1$, $2$, $2$, $3$, $3$, $4$, $4$, $5$ and $6$, and it is
the smallest group supporting an icosahedral representation (since
$A_5$ has no $2$-dimensional representations), meaning that our
criterion unfortunately does not apply to checking the icosahedral
case.  In fact, one knows Artin's conjecture for {\em all} induced
representations of this group; while they are not all monomial, the
only exceptions come from a pair of tetrahedral representations, for
which we have the Langlands-Tunnell theorem.  Even with this added
information, we cannot rule out the possibility of a simple pole with
undetectably small residue at a zero of the $L$-function of the
$6$-dimensional representation.  More precisely, we find with GAP that
the induced representations are spanned by the twelve vectors
\begin{equation}
\begin{aligned}
&(1,0,0,0,0,0,0,0,0)\qquad
&(0,0,0,0,0,0,0,1,0)\qquad
&(1,0,0,0,0,1,0,0,0)\\
&(0,0,0,1,1,0,0,0,0)\qquad
&(0,0,0,1,0,1,0,1,0)\qquad
&(0,0,0,0,1,1,0,1,0)\\
&(0,0,0,1,1,1,0,0,0)\qquad
&(0,0,0,0,0,0,0,0,1)\qquad
&(0,0,0,0,0,0,1,0,1)\\
&(0,1,0,0,0,0,1,0,1)\qquad
&(0,0,1,0,0,0,1,0,1)\qquad
&(0,1,1,0,0,0,0,0,1).
\end{aligned}
\end{equation}
The first seven of these are the monomial representations lifted from
$\SL_2(\F_5)/\{\pm I\}\cong A_5$, while the others give ``new''
information.  One easily checks that for $\rho$ the $6$-dimensional
representation, \eqref{eq:monohyp} fails with $\chi_1$ corresponding to
any of the vectors $(0,-1,0,0,0,0,0,0,1)$, $(0,0,-1,0,0,0,0,0,1)$
and $(0,0,0,0,0,0,-1,0,1)$, i.e.\ the representations of dimension $2$
and one of dimension $4$ can hide a pole at a zero of $L(s,\rho)$.
This shows in a strong sense that information from induced
representations is in general insufficient to show Artin's conjecture.

However, all is not lost concerning icosahedral representations.
For a given icosahedral $\rho$, the {\em adjoint square} $\Ad(\rho)$
is a $3$-dimensional representation with image isomorphic to $A_5$.
A result of Flicker \cite{flicker} implies that modularity of
$\rho$ is equivalent to that of $\Ad(\rho)$.  (In fact, modularity
of all representations of the underlying group follows from that of
$\Ad(\rho)$ and its Galois conjugate, by known cases of functoriality;
see \cite{wang}.)  Combining this fact with the $\GL(3)$ converse
theorem, one could give a converse theorem for $\GL(2)$ using analytic
properties of $L(s,\Ad(\rho)\otimes\chi)$ for Dirichlet characters $\chi$.
Weissman, in his undergraduate thesis \cite{weissman}, used this idea
to give indirect evidence for the modularity of an even icosahedral
representation.  By Propositions \ref{prop:monofunc} and \ref{prop:monoex}
we see that in principle we may directly verify the holomorphy of these
$L$-functions up to finite height.  Moreover, an ``effective'' version of
the $\GL(3)$ converse theorem (requiring, say, meromorphy of all twists
and holomorphy of a finite number in a bounded region) would suffice to
give an algorithm for verifying the conjecture in the icosahedral case.

Unfortunately, there is the more practical problem that totally real
$A_5$ fields (those that give rise to even icosahedral representations)
are very rare; the smallest known discriminant is far too large to test
with current computers.  Thus, for the $A_5$ examples that we consider in
Section \ref{numerics}, the Artin conjecture is already known.  To test
our criterion, we consider instead some examples of $S_5$ extensions,
which exist in much greater abundance.

Finally, we note that in the course of verifying Artin's conjecture,
the information that we collect implies that the zeros of each
$L(s,\rho)$ are simple and lie on the line $\Re(s)=\frac12$.  Thus, in
the process we also verify the Riemann hypothesis for $\zeta_K(s)$.
Interestingly, we do not need to establish the holomorphy of all
$L(s,\rho)$ in order to do this; it is enough, for example, that they
have at most simple poles.  More precisely, in order to check the
Riemann hypothesis around a generic zero of $\zeta_K(s)$ we need to have
\begin{equation}
\Tr\rho\ne\chi_1+2\chi_2\mbox{ for virtual characters }
\chi_i\ne 0\mbox{ with }\langle\chi_i,\sigma\rangle\ge 0
\mbox{ for all monomial }\sigma,
\end{equation}
which is a weaker condition than almost monomiality.  In particular, we
may still check the Riemann hypothesis for $\SL_2(\F_5)$ extensions.

\section{Locating zeros via the explicit formula}
\label{sec:explicit}
Let notation be as in the introduction, and define numbers $c_n$ by
$-\frac{L'}{L}(s) = \sum_{n=1}^{\infty}c_nn^{-s}$,
i.e.\ $c_n=(\log p)\sum_{j=1}^r\alpha_{p,j}^k$ for $n=p^k$ a prime power,
and $c_n=0$ otherwise.  Further, we enumerate the zeros of $\Lambda(s)$ as
$\rho_n=\frac12+i\gamma_n$ for $n\in\Z$, repeated with multiplicity.
Weil's explicit formula relates the sequences $\{c_n\}$ and $\{\gamma_n\}$.
Precisely, let $g\in C^1_c(\R)$
be a differentiable function of compact support such that its Fourier
transform $h(t):=\int_{-\infty}^{\infty}g(x)e^{ixt}\,dx$ is real for
$t\in\R$.  Then
\begin{equation}
\begin{aligned}
\sum_{n\in\Z}&h(\gamma_n)
-2\Re\sum_{k=1}^m h\!\left(-i\!\left(\frac12+\lambda_k\right)\!\right)\\
&=g(0)\log N
+2\Re\!\left[\sum_{j=1}^r\frac1{2\pi}\int_{-\infty}^{\infty}
\frac{\Gamma_{\R}'}{\Gamma_{\R}}\!\left(\frac12+\mu_j+it\right)\!h(t)\,dt
-\sum_{n=1}^{\infty}\frac{c_n}{\sqrt n}g(\log n)\right],
\end{aligned}
\end{equation}
This follows from the Cauchy integral formula and the functional
equation; see \cite{rudnick-sarnak}.
Note that all terms of the formula may be put in terms of $g$; in
particular,
\begin{equation}
\frac1{2\pi}\int_{-\infty}^{\infty}
\frac{\Gamma_{\R}'}{\Gamma_{\R}}\!\left(\frac12+\mu+it\right)\!h(t)\,dt
=\frac12\int_0^{\infty}\log\bigl(\pi e^{\gamma}(e^{2x}-1)\bigr)
\,d\bigl(g(x)e^{-(1/2+\mu)x}\bigr).
\end{equation}
This form is convenient for computation, since $g$ has compact support.

The important thing to note is that given a list of the $c_n$ for
$n\le e^X$, the explicit formula gives a method for evaluating $\sum_n
h(\gamma_n)$ for essentially any function $h$ whose Fourier transform
is supported in $[-X,X]$.  When $X$ is large, we may choose $h$ to be
narrowly concentrated around any particular point, and thus resolve
features of the spectrum in places where the density of zeros is not too
large compared to $X$; a variant of this technique, with explicit test
functions (not of compact support), was worked out by Omar \cite{omar}
to estimate the lowest zero of some Dedekind zeta functions.  For a
fixed support $[-X,X]$, there is a canonical way of choosing a ``best''
test function, by a method developed for the Selberg trace formula in
\cite{bs}.  In order to use the method, which depends crucially on a
positivity argument, it is necessary to assume the Riemann hypothesis
for our given $L$-function.  With that caveat, we recall briefly the
construction from \cite{bs}.

For $t_0\in\R$, let ${\mathcal C}(X,t_0)$ be the
class of functions $h$ as above, with the corresponding $g$ supported
in $[-X,X]$, and the additional restrictions $h(t)\ge 0$ for
$t\in\R$ and $h(t_0)=1$.  Define
\begin{equation}
F_X(t_0) := \inf_{h\in{\mathcal C}(X,t_0)}\sum_{n\in\Z}h(\gamma_n).
\label{eq:FXdef}
\end{equation}
Then as $X\to\infty$, $F_X$ tends pointwise to the characteristic
function of the zeros.  Moreover, if $F_X(t_0)<1$ for any value of $X$
then $t_0$ cannot be the ordinate of a zero.  Thus, by evaluating $F_X$
we can find provable intervals in which the zeros must lie.

Although the definition of $F_X$ is abstract, it is easy to construct
concrete families of functions that closely approximate any desired
function.  For instance, let $M$ be a large integer, $\delta=X/2M$ and set
\begin{equation}
h(t) = \left(\frac{\sin\delta t/2}{\delta t/2}\right)^4
\left(a_0+\sum_{n=1}^{M-1}
\bigl(a_n\cos\delta nt+b_n\sin\delta nt\bigr)\right)^2,
\end{equation}
for arbitrary real numbers $a_n,b_n$.
(For self-dual $L$-functions, we restrict to even test
functions, i.e.\ all $b_n=0$, and divide the final formula by $2$.)
On the other side of the Fourier transform, this corresponds to taking
$g=f*f$, where $f$ linearly interpolates arbitrary values at multiples
of $\delta$.

The sum over zeros in \eqref{eq:FXdef} is then a positive definite
quadratic form in the numbers $a_n$ and $b_n$.  To compute the matrix
of the form essentially involves computing the explicit formula for
functions $g$ that are translates of a fixed function of small compact
support.  That requires almost no extra work, since we may compute the
formula for all localized test functions simultaneously.  Once the
matrix is known, the infemum in \eqref{eq:FXdef} over this restricted
class of test functions is easily found as the minimum of the quadratic
form subject to the linear constraint $h(t_0)=1$.  This involves
inverting the matrix, after which the minimum may be found quickly
for many different values of $t_0$.

For an $L$-function of degree $r$ and conductor $N$, the
density of zeros at height $T$ is roughly $\frac1{2\pi}\log
N\!\left(\frac{T}{2\pi}\right)^r$.  Therefore, in order to resolve features
around height $T$, the uncertainty principle says we should know the
numbers $c_n$ for $n$ up to about $N\!\left(\frac{T}{2\pi}\right)^r$.
In the self-dual case, the extra division by $2$ replaces this by its
square root; thus, the complexity is on par with that of the approximate
functional equation or the algorithm of Section \ref{sec:rigorous},
although it is much more sensitive to the local spacing of zeros.
(Heuristic arguments based on experiments and random matrix theory
\cite{odlyzko} indicate that the minimum gap between zeros can be
arbitrarily small relative to the mean value; although such small gaps
are expected to be very rare, we could in principle need many more
coefficients than for the ``typical'' zero at height $T$.)  In practice,
the explicit formula is clean and easy to implement since there are
no error terms to estimate with functions of compact support.  It is
particularly well-suited to finding low zeros or to situations where
the numbers $c_n$ may be computed quickly, as is the case for Artin
$L$-functions; cf.\ Section \ref{numerics:coeff}.

As mentioned above, the minimization procedure requires assuming the
Riemann hypothesis.  If one is willing to do so, the method may be made
completely rigorous, and may even be used in place of Turing's method
for verifying Artin's conjecture.  However, it is more natural to use
it as a quick check in order to fine tune and validate the subsequent
rigorous methods.  In fact, it is helpful to assume Artin's conjecture
and apply the method to the irreducible Artin $L$-functions directly.
That thins out the spectrum, making it easier to isolate individual
zeros.  We have carried out this procedure for a few examples in Section
\ref{numerics:explicit}.

\section{Turing's method} \label{turing}
Turing's method for verifying the Riemann hypothesis is described
well in his paper \cite{turing}, although there are some errors
in the details that were later corrected by Lehman \cite{lehman}.
The method has subsequently been extended to Dirichlet $L$-functions
by Rumely \cite{rumely} and Dedekind zeta functions by Tollis
\cite{tollis}\footnote{Tollis applied his method to cubic and quartic
fields.  In these cases, there is a slight advantage in passing to the
normal closure and separating into irreducible Artin $L$-functions,
as we have done for the $A_5$ cases in Section \ref{numerics}.}.
Our contribution is to work out the details necessary to apply it to an
arbitrary $L$-function with simple zeros.

Our argument essentially follows that of Turing.  We begin by setting
some notation to be used only in this section.  For $t$ not the
ordinate of a zero or pole of $\Lambda$, let
\begin{equation}
S(t):=\frac1{\pi}\Im\int_{\infty}^{1/2}\frac{L'}{L}(\sigma+it)\,d\sigma.
\end{equation}
By convention, we make $S(t)$ upper semi-continuous, i.e.\ when $t$ is
the ordinate of zero or pole, we define
$S(t)=\lim_{\varepsilon\to0^+}S(t+\varepsilon)$.

Next, for $t_1<t_2$ let $N(t_1,t_2)$ denote the net number of zeros
with imaginary part in $(t_1,t_2]$,
counting multiplicity.
When neither $t_1$ nor $t_2$ is the ordinate of a zero or pole, we may
calculate $N(t_1,t_2)$ using the argument principle, as in \eqref{eq:argp1}.
Let $C$ be the rectangle with corners at
$2+it_1$, $2+it_2$, $-1+it_2$, $-1+it_1$, with
counter-clockwise orientation, $H$ the half plane
$\bigl\{s\in\C:\Re(s)\ge\frac12\bigr\}$, and $H^c$ its complement.
Note that by the functional equation, we have
$\overline{\frac{\Lambda'}{\Lambda}(s)}=
-\frac{\Lambda'}{\Lambda}(1-\bar{s})$.  Hence,
\begin{equation}
\begin{aligned}
N(t_1,t_2)
&=\frac1{2\pi}\Im\int_C\frac{\Lambda'}{\Lambda}(s)\,ds
=\frac1{2\pi}\Im\!\left(
\int_{C\cap H}\frac{\Lambda'}{\Lambda}(s)\,ds
-\overline{\int_{C\cap H^c}\frac{\Lambda'}{\Lambda}(s)\,ds}\right)\\
&=\frac1{2\pi}\Im\!\left(
\int_{C\cap H}\frac{\Lambda'}{\Lambda}(s)\,ds
+\int_{C\cap H^c}\frac{\Lambda'}
{\Lambda}(1-\bar{s})\,d\bar{s}\right)\\
&=\frac1{\pi}\Im\int_{C\cap H}\frac{\Lambda'}{\Lambda}(s)\,ds
=\frac1{\pi}\Im\int_{C\cap H}\frac{\gamma'}{\gamma}(s)\,ds
+\frac1{\pi}\Im\int_{C\cap H}\frac{L'}{L}(s)\,ds.
\label{Ntargprinc}
\end{aligned}
\end{equation}
Now for the integral of $L'/L$ we move the right edge of the contour out
to $\infty$, where the integrand vanishes.  We thus obtain
\begin{equation}
N(t_1,t_2) = \frac1{\pi}\Im\log\gamma(s)\bigr|_{\frac12+it_1}^{\frac12+it_2}
+S(t_2)-S(t_1).
\end{equation}

We select a particular branch of $\log\gamma\bigl(s)$
by using the principal branch of $\log\Gamma$.  With this choice, set
\begin{equation}
\Phi(t) := \frac1{\pi}\left[
\arg\epsilon+\frac{\log N}2t
-\frac{\log\pi}2\!\left(rt+\Im\sum_{j=1}^r\mu_j\right)
+\Im\sum_{j=1}^r
\log\Gamma\!\left(\frac{1/2+it+\mu_j}2\right)\right]
\label{gammabranch}
\end{equation}
and
\begin{equation}
N(t) := \Phi(t)+S(t).
\end{equation}
Then $N(t_1,t_2)=N(t_2)-N(t_1)$.
Note that if $L$ is self-dual and vanishes to order $\le 1$ at $\frac12$ then
$N(t)=N(0,t)$.  In the general case, although we still have
$N(t)\in\Z$, there is no standard reference point, so only changes in
$N(t)$ are meaningful.  (Put another way, the branch of $\log\gamma$
chosen in \eqref{gammabranch} is non-canonical.)  For large $t$, $\Phi(t)$
may be evaluated quickly by an effective version of Stirling's formula:
\begin{equation}
\Im\log\Gamma(z)=\Im\!\left[\left(z-\frac12\right)\log\frac{z}{e}\right]
+\Theta\!\left(\frac1{8|\Im(z)|}\right)
\qquad\mbox{for }z\in\C\setminus\R,
\end{equation}
where the notation $f=\Theta(g)$ means $|f|\le g$.

Turing's method is as follows.  Recall that
$\Lambda\bigl(\frac12+it\bigr)$ is real valued.  Thus, if we have an
accurate procedure to compute $\Lambda(s)$ then we may locate all
simple zeros on the line $\Re(s)=\frac12$ by observing its sign
changes.  If it turns out that all of the zeros between ordinates $t_1$
and $t_2$ are simple and on the line, then we can deduce the Riemann
hypothesis in that interval by computing $N(t_1,t_2)$ (minus the
contribution from any poles between $t_1$ and $t_2$) and finding the
same number of sign changes over the interval.

To compute $N(t_1,t_2)$, we could evaluate \eqref{Ntargprinc}
numerically.  However, this would require many evaluations of
$\Lambda(s)$ and would be difficult to carry out rigorously.  Fortunately,
Turing devised a simpler method, based on the fact (first due to Littlewood
for $\zeta(s)$) that $S(t)$ has mean value $0$.  Thus, the graph of
$N(t_0,t)-\Phi(t)$ for any fixed $t_0$ oscillates around a constant value;
if we were to plot the same function using the {\em measured} number of
zeros in $(t_0,t]$, then any zeros that we had missed would be obvious
as jumps in the graph.

This can be made rigorous as follows.  Let $t_0$ be a large number that
is not the ordinate of a zero or pole, and assume that between
ordinates $t_0-h$ and $t_0+h$ (for some $h>0$), we have located several
zeros of $\Lambda(s)$, i.e.\ we have found small intervals $(a_n,b_n)$
such that $\Lambda\bigl(\frac12+ia_n\bigr)$ and
$\Lambda\bigl(\frac12+ib_n\bigr)$ have opposite sign.
Let $\Nlt(t_0,t)$ (resp.\ $\Nrt(t_0,t)$) be the step function which is
upper semi-continuous, increases by $1$ at each $a_n$ (resp.\ $b_n$)
and vanishes at $t=t_0$.  We then have
\begin{equation}
N(t)\le N(t_0)+\Nlt(t_0,t)\mbox{ for }t\le t_0
\qquad\mbox{and}\qquad
N(t)\ge N(t_0)+\Nrt(t_0,t)\mbox{ for }t\ge t_0.
\label{Ntbounds}
\end{equation}
From these, we can deduce upper and lower bounds for $N(t_0)$;
integrating \eqref{Ntbounds}, we get
\begin{equation}
N(t_0)h+\int_{t_0}^{t_0+h}\Nrt(t_0,t)\,dt
\le \int_{t_0}^{t_0+h}N(t)\,dt
=\int_{t_0}^{t_0+h}\Phi(t)\,dt
+\int_{t_0}^{t_0+h}S(t)\,dt
\label{Nupper}
\end{equation}
and
\begin{equation}
N(t_0)h+\int_{t_0-h}^{t_0}\Nlt(t_0,t)\,dt
\ge \int_{t_0-h}^{t_0}N(t)\,dt
=\int_{t_0-h}^{t_0}\Phi(t)\,dt
+\int_{t_0-h}^{t_0}S(t)\,dt.
\label{Nlower}
\end{equation}

If we have in fact located all zeros in the interval $(t_0-h,t_0+h)$
with some amount of precision (as measured by the size of the intervals
$(a_n,b_n)$), then we can expect these bounds to be
close to the truth.  Moreover, if we have effective upper and lower
bounds for the integral of $S(t)$, then for $h$ large enough,
\eqref{Nupper} and \eqref{Nlower} will bound a single integer,
i.e.\ we can unambiguously determine $N(t_0)$.  Doing this for two
different values $t_0=t_1, t_2$, we obtain $N(t_1,t_2)$.

One nice feature of Turing's method is that precise knowledge of the
zeros is only required in the short intervals around $t_1$ and $t_2$,
and even there one can make a trade-off between the precision of the
zeros and the length $h$ of the interval.  For the bulk of the zeros
between $t_1$ and $t_2$ it suffices to observe the sign changes.

The remainder of this section is devoted to bounding $\int S(t)\,dt$,
cf.\ Theorem \ref{Stbound} below.  Our starting point is the following
formula, obtained by Littlewood's box principle (see
\cite[\S9.9]{titchmarsh}):
\begin{equation}
\pi\int_{t_1}^{t_2}S(t)\,dt
=\int_{1/2}^{\infty}\log|L(\sigma+it_2)|\,d\sigma
-\int_{1/2}^{\infty}\log|L(\sigma+it_1)|\,d\sigma.
\label{Stformula}
\end{equation}

\begin{lemma}
\label{Lupperbound}
Let notation be as above, and set
$B:=\sup_{\Re(s)=3/2}|L(s)|^2$.
Then, for $s$ in the strip $\bigl\{s\in\C:-\frac12\le \Re(s)\le \frac32\bigr\}$,
\begin{equation}
|L(s)|^2\le B|\chi(s)Q(s)|\left|\frac{P(s+1)^2P(s-2)}{P(s)^2P(s-1)}\right|.
\end{equation}
\end{lemma}
\begin{remark}
The power of $|Q(s)|$ in the above is not optimal; for
$\Re(s)=\frac12$, the ``convexity bound'' says that we can put instead
$|Q(s)|^{1/2+\varepsilon}$, with a constant depending on $\varepsilon$
(see \cite{iwaniec-sarnak}), while the Lindel\"of hypothesis would have
$|Q(s)|^{\varepsilon}$.  Our present choice permits us to avoid
Stirling's formula in the proof, and thus obtain a clean bound that
is uniform in all parameters.
\end{remark}
\begin{proof}
We consider first the case when $L(s)$ is entire.
Set
\begin{equation}
F(s):=L(s)\Lbar(1-s) = \chi(s)^{-1}L(s)^2.
\end{equation}
Plugging in the definition of $\chi(s)$,
\begin{equation}
|F(\sigma+it)|=|L(\sigma+it)|^2
\left|\frac{\gamma(\sigma+it)}{\gammabar(1-\sigma-it)}\right|
=|L(\sigma+it)|^2
\left|\frac{\gamma(\sigma+it)}{\gamma(1-\sigma+it)}\right|.
\end{equation}
Note that when $\sigma=\frac12+\mbox{a positive integer}$, the ratio of
$\gamma$ factors reduces to a polynomial; in particular,
\begin{equation}
|F(3/2+it)|=|L(3/2+it)|^2|Q(-1/2+it)| \le B|Q(3/2+it)|.
\label{F32bound}
\end{equation}
The inequality holds since $\Re(\mu_j)\ge -\frac12$ for all $j$.
Next, from the functional equation we have $F(s)=\Fbar(1-s)$, so that
\begin{equation}
|F(\sigma+it)|=|F(1-\sigma+it)|.
\end{equation}
Hence, by \eqref{F32bound},
\begin{equation}
|F(-1/2+it)|\le B|Q(-1/2+it)|.
\end{equation}
Thus, the function $F(s)/Q(s)$ is bounded by $B$ on the lines
$\Re(s)=-\frac12$ and $\Re(s)=\frac32$.  Note that although $Q(s)$ has zeros,
$F(s)$ has trivial zeros at the same points; in fact
\begin{equation}
\frac{F(s)}{Q(s)}=\frac{\Lambda(s)\Lbar(1-s)}{\gamma(s)Q(s)}
=\frac{\Lambda(s)\Lbar(1-s)}{\gamma(s+2)}.
\label{FPholo}
\end{equation}
Since $F$ has finite order, it follows from the Phragmen-Lindel\"of
theorem that $|F(s)|\le B|Q(s)|$ for all $s$ in the strip.

If $L(s)$ has poles then the above argument breaks down since
$F(s)/Q(s)$ is not holomorphic in the strip.  In fact, for each $k$ we
get three poles, one at $1+\lambda_k$ and two at $\lambda_k$, as \eqref{FPholo}
shows.  To compensate for this,
we consider $F(s)P(s)^2P(s-1)$ in the above, in place of $F(s)$.
One checks that $|s^2(s-1)|\le |(s+1)^2(s-2)|$ on the lines
$\Re(s)=-1/2$ and $\Re(s)=3/2$, so that
\begin{equation}
|F(s)P(s)^2P(s-1)|\le B|Q(s)P(s)^2P(s-1)| \le B|Q(s)P(s+1)^2P(s-2)|.
\end{equation}
Further, the ratio $\frac{F(s)P(s)^2P(s-1)}{Q(s)P(s+1)^2P(s-2)}$
is holomorphic in the strip, so we may proceed as above.  The lemma follows.
\end{proof}

\begin{lemma}
\label{technical}
Suppose that
\begin{equation}
(t+\Im(\mu_j))^2\ge (5/2+\Re(\mu_j))^2+X^2
\quad\mbox{for some }X>5\mbox{ and all }j=1,\ldots,r.
\label{tmuhyp}
\end{equation}
Then
\begin{enumerate}
\item[1.]
For $\sigma\in [1/2,5/2]$,
\begin{equation}
-r\left(\frac1{2\sqrt2X}+\frac{4/\pi^2+1/4}{X^2}\right) \le
\Re\frac{\gamma'}{\gamma}(\sigma+it)
-\frac12\log\left|Q\!\left(\frac32+it\right)\right|
\le \frac{4r}{\pi^2X^2}.
\label{gammabound}
\end{equation}
\item[2.]
For $\sigma\in [-1/2,3/2]$,
\begin{equation}
\Re\frac{Q'}{Q}(\sigma+it)\le \frac{r}{\sqrt2 X}.
\end{equation}
\item[3.]
For all $\sigma$,
\begin{equation}
\Re\frac{P'}{P}(\sigma+it)\le\frac{\max(\sigma m,0)}{X^2}.
\end{equation}
\end{enumerate}
\end{lemma}
\begin{proof}
1. We have
\begin{equation}
\frac{\gamma'}{\gamma}(\sigma+it)
=\frac12\log\frac{N}{\pi^r}+\frac12\sum_{j=1}^r
\frac{\Gamma'}{\Gamma}\!\left(\frac{\sigma+it+\mu_j}2\right).
\end{equation}
We apply the Stirling-type estimate \cite{lehman}
\begin{equation}
\frac{\Gamma'}{\Gamma}(z) = \log z-\frac1{2z}
+\Theta\!\left(\frac{2/\pi^2}{|\Im(z)^2-\Re(z)^2|}\right)
\qquad\mbox{for }\Re(z)\ge 0.
\label{Gammalogderiv}
\end{equation}
This yields
\begin{equation}
\begin{aligned}
&\Re\frac{\gamma'}{\gamma}(\sigma+it)
-\frac12\log\left|Q\!\left(\frac32+it\right)\right|\\
&=-\frac12\sum_{j=1}^r\left[
\log\left|\frac{3/2+it+\mu_j}{\sigma+it+\mu_j}\right|
+\Re\frac1{\sigma+it+\mu_j}
+\Theta\!\left(\frac{8/\pi^2}{|(t+\Im(\mu_j))^2-(\sigma+\Re(\mu_j))^2|}\right)
\right].
\end{aligned}
\end{equation}

For the lower bound in \eqref{gammabound} we need an upper bound for
the expression in brackets.  By hypothesis, the $\Theta$ term is
bounded by $\frac8{\pi^2 X^2}$.  For the others, put
$\sigma+it+\mu_j=x+iy$, $\beta=3/2-\sigma$, so that $|\beta|\le 1$ and
$x$ and $y$ are constrained by $y^2\ge x^2+X^2$, $x\ge 0$.  Then, using
the inequality $\log(1+u)\le u$, we have
\begin{equation}
\begin{aligned}
\log\left|\frac{x+\beta+iy}{x+iy}\right|+\Re\frac1{x+iy}
&=\frac12\log\!\left(1+\frac{2\beta x+\beta^2}{x^2+y^2}\right)+
\frac{x}{x^2+y^2}\\
&\le \frac{(\beta+1)x+\beta^2/2}{x^2+y^2}
\le\frac{2x+1/2}{2x^2+X^2} \le \frac1{\sqrt{2}X}+\frac1{2X^2}.
\end{aligned}
\end{equation}
There are $r$ such terms, and the lower bound follows after multiplying
by $-\frac12$.

The upper bound is similar, but uses the second order inequality
$\log(1+u)\ge u-u^2/2$.  We omit the details.

2. Similarly,
\begin{equation}
\Re\frac{Q'}{Q}(\sigma+it)
=\sum_{j=1}^r\Re\frac1{\sigma+\mu_j+it}
\le\sum_{j=1}^r\frac{\sigma+\Re(\mu_j)}
{2(\sigma+\Re(\mu_j))^2+X^2}\le \frac{r}{\sqrt2X}.
\end{equation}

3.
\begin{equation}
\Re\frac{P'}{P}(\sigma+it)
=\sum_{k=1}^m\Re\frac1{\sigma+it-\lambda_k}
\le\sum_{k=1}^m\frac{\max(\sigma,0)}
{\sigma^2+X^2}\le \frac{\max(\sigma m,0)}{X^2}.
\end{equation}
(Note that this bound is of faster decay than estimates 1 and 2.
That is because we have control over the real parts of the poles, while
nothing prevents $\Re(\mu_j)$ from being of comparable size to $X$.
One could also obtain an $O\bigl(X^{-2}\bigr)$ bound in 1 and 2, with
a constant depending on the $\mu_j$.)
\end{proof}

\begin{lemma}
\label{hardlemma}
Let $w$ be a complex number with $|\Re(w)|\le \frac12$.  Then
\begin{equation}
\int_0^1\log\left|
\frac{(x+1+w)(x+1-\overline{w})}{(x+w)(x-\overline{w})}\right|dx
\le (\log 4)\Re\!\left(\frac1{1+w}+\frac1{1-\overline{w}}\right).
\label{hardineq}
\end{equation}
\end{lemma}
\begin{proof}[Proof (sketch)]
Note first that equality is attained at $w=0$.
Using the principal branch of the logarithm, set
\begin{equation}
f(w)=\int_0^1\log\!\left(
\frac{(x+1+w)(x+1-w)}{(x+w)(x-w)}\right)dx
\quad\mbox{and}\quad
g(w)=\frac1{1+w}+\frac1{1-w}.
\end{equation}
These define analytic functions on $\C\setminus\R$, with real parts
extending continuously to $\R$.  Further, \eqref{hardineq}
is equivalent to the assertion
\begin{equation}
\Re(f(w)-(\log 4)g(w)) \le 0 \quad\mbox{for }|\Re(w)|\le\frac12.
\label{equivassertion}
\end{equation}

Note that $f(w)$ and $g(w)$ are each asymptotic to $-2/w^2$ as
$|w|\to\infty$.  Thus, \eqref{equivassertion} holds for $\Im(w)$
sufficiently large; we check that in fact $|\Im(w)|\ge 2$ is enough.
By symmetry and the maximum modulus principle applied to the function
$e^{f(w)-(\log 4)g(w)}$ on the rectangle with corners at $\pm\frac12$
and $\pm\frac12+2i$, it suffices to check \eqref{equivassertion} on the
real axis and for $\Re(w)=\frac12$.  On the real axis we calculate the
integral explicitly and verify the inequality using calculus.  For
$\Re(w)=\frac12$ the inequality is strict, so we may verify it
computationally for $0\le \Im(w) \le 2$.
\end{proof}

\begin{lemma}
\label{zetalemma}
For $\sigma > \theta+1$, define
\begin{equation}
z_{\theta}(\sigma):=
\left(\frac{\zeta(2\sigma+2\theta)\zeta(2\sigma-2\theta)}
{\zeta(\sigma+\theta)\zeta(\sigma-\theta)}\right)^{1/2}
\quad\mbox{and}\quad
Z_{\theta}(\sigma):=
\bigl(\zeta(\sigma+\theta)\zeta(\sigma-\theta)\bigr)^{1/2}.
\end{equation}
Then
\begin{equation}
r\frac{z_{\theta}'}{z_{\theta}}(\sigma)\ge
\Re\frac{L'}{L}(\sigma+it)\ge
r\frac{Z_{\theta}'}{Z_{\theta}}(\sigma)
\label{zetapart1}
\end{equation}
and
\begin{equation}
z_{\theta}(\sigma)^r\le |L(\sigma+it)|\le Z_{\theta}(\sigma)^r.
\label{zetapart2}
\end{equation}
\end{lemma}
\begin{proof}
From \eqref{eulerprod} we have
\begin{equation}
r\frac{Z_{\theta}'}{Z_{\theta}}(\sigma)-\Re\frac{L'}{L}(\sigma+it)
=\sum_p\sum_{k=1}^{\infty}p^{-k\sigma}\log{p}
\left(\Re\sum_{j=1}^r\bigl(\alpha_{p,j}p^{-it}\bigr)^k
-\frac r2\bigl(p^{k\theta}+p^{-k\theta}\bigr)\right).
\end{equation}
Pairing the terms for $\alpha_{p,j}$ and $\alpha_{p,j}'$, we see that each
summand is $\le 0$, from which the second inequality of \eqref{zetapart1}
follows.  The first inequality is similar.  For \eqref{zetapart2},
integrate \eqref{zetapart1} from $\sigma$ to $\infty$.
\end{proof}

\begin{theorem}
\label{Stbound}
Suppose $t_1$ and $t_2$ satisfy \eqref{tmuhyp}, and set
\begin{equation}
c_{\theta}:=\log Z_{\theta}\!\left(\frac32\right)
+\int_{3/2}^{\infty}
\log\frac{Z_{\theta}(\sigma)}{z_{\theta}(\sigma)}\,d\sigma
-\int_{3/2}^{5/2}\log z_{\theta}(\sigma)\,d\sigma
+(\log 4)\frac{z_{\theta}'}{z_{\theta}}\!\left(\frac32\right),
\end{equation}
where $z_{\theta}$ and $Z_{\theta}$ are as in Lemma \ref{zetalemma}.
(In particular, $c_0 \lessapprox 5.65055$.)  Then
\begin{equation}
\pi\int_{t_1}^{t_2}S(t)\,dt \le
\frac14\log\left|Q\!\left(\frac32+it_2\right)\right|
+\left(\log2-\frac12\right)
\log\left|Q\!\left(\frac32+it_1\right)\right|
+c_{\theta}r+\frac{r}{\sqrt2(X-5)}.
\end{equation}
\end{theorem}
\begin{remark}
Note that there is no assumption on the order of $t_1$ and $t_2$, so one
obtains a lower bound as well by reversing their roles.
\end{remark}
\begin{proof}
By \eqref{Stformula}, we need upper and lower bounds for
$\int_{1/2}^{\infty}\log|L(\sigma+it)|\,d\sigma$.
For the upper bound, we use Lemma \ref{Lupperbound}:
\begin{equation}
\begin{aligned}
\int_{1/2}^{\infty}\log|L(\sigma+it)|&\,d\sigma
=\int_{1/2}^{3/2}\log|L(\sigma+it)|\,d\sigma
+\int_{3/2}^{\infty}\log|L(\sigma+it)|\,d\sigma\\
&\le\frac12\log B
+\frac12\int_{1/2}^{3/2}\log\left|
\frac{\gamma(1-\sigma+it)Q(\sigma+it)}{\gamma(\sigma+it)}
\right|d\sigma\\
&+\frac12\int_{1/2}^{3/2}\log\left|
\frac{P(\sigma+1+it)^2P(\sigma-2+it)}{P(\sigma+it)^2P(\sigma-1+it)}
\right|d\sigma
+\int_{3/2}^{\infty}\log|L(\sigma+it)|\,d\sigma.
\end{aligned}
\end{equation}
We bound the first and last terms with Lemma \ref{zetalemma}:
\begin{equation}
\sup_{\Re(s)=3/2}|L(s)|
+\int_{3/2}^{\infty}\log|L(\sigma+it)|\,d\sigma
\le \log Z_{\theta}\!\left(\frac32\right)
+\int_{3/2}^{\infty}\log Z_{\theta}(\sigma)\,d\sigma.
\end{equation}
For the second term, we replace $\sigma$ by $2-\sigma$ in the top
$\gamma$ factor, and use the recurrence for $\gamma$ to get
\begin{equation}
\frac12\int_{1/2}^{3/2}\log\left|
\frac{\gamma(\sigma+1+it)}{\gamma(\sigma+it)}
\frac{Q(\sigma+it)}{Q(\sigma-1+it)}\right|d\sigma.
\label{uppersecond}
\end{equation}
By the mean value theorem, the integral equals
$\Re\bigl[\frac{\gamma'}{\gamma}(\sigma^*+it)
+\frac{Q'}{Q}(\sigma^*-1+it)\bigr]$
for some $\sigma^* \in [1/2,5/2]$.  Thus, by Lemma \ref{technical},
\eqref{uppersecond} is at most
\begin{equation}
\frac14\log\left|Q\!\left(\frac32+it\right)\right|
+r\left(\frac2{\pi^2X^2}+\frac1{2\sqrt 2X}\right).
\end{equation}
Similarly, we see that the third term is bounded by $\frac{5m}{2X^2}$.

We now turn to the lower bound.  This part is more delicate
since we must take into account the contribution of zeros near
$\frac12+it$.  We use Turing's idea of comparing $\log|L(s)|$ to
$\log|L(s+1)|$; the difference between these looks like a value of the
logarithmic derivative, which we can make precise with the help of
Lemma \ref{hardlemma}.

Proceeding, we first clear the poles of $\Lambda$ by writing
$F(s)=\Lambda(s)P(s)P(s-1)$.  This function then has the
Weierstrass-Hadamard product
\begin{equation}
F(s) = e^{as+b}\prod_{\rho}\left(1-\frac{s}{\rho}\right)
e^{s/\rho},
\label{hadamard}
\end{equation}
where $\rho$ runs over the zeros of $\Lambda$ and
$\Re(a) = -\sum_{\rho}\Re\bigl(\frac1{\rho}\bigr)$.
Next, we split the integral as follows:
\begin{equation}
\begin{aligned}
\int_{1/2}^{\infty}&\log|L(\sigma+it)|\,d\sigma
=\int_{1/2}^{3/2}\log\left|\frac{F(\sigma+it)}{F(\sigma+1+it)}\right|d\sigma
+\int_{1/2}^{3/2}\log\left|\frac{\gamma(\sigma+1+it)}
{\gamma(\sigma+it)}\right|d\sigma\\
&+\int_{1/2}^{3/2}\log\left|\frac{P(\sigma+1+it)}
{P(\sigma-1+it)}\right|d\sigma
+ \int_{3/2}^{5/2}\log|L(\sigma+it)|\,d\sigma
+\int_{3/2}^{\infty}\log|L(\sigma+it)|\,d\sigma.
\end{aligned}
\end{equation}
The second term may be estimated, as above, by Lemma
\ref{technical} and the mean value theorem:
\begin{equation}
\int_{1/2}^{3/2}\log\left|\frac{\gamma(\sigma+1+it)}
{\gamma(\sigma+it)}\right|d\sigma
\ge \frac12\log\left|Q\!\left(\frac32+it\right)\right|
-r\left(\frac1{2\sqrt2X}+\frac{4/\pi^2+1/4}{X^2}\right).
\end{equation}
The third term is positive since
$\bigl|\frac{P(\sigma+1+it)}{P(\sigma-1+it)}\bigr|\ge 1$ for
$\sigma\ge 0$.
The fourth and fifth terms are handled by Lemma \ref{zetalemma}.

As for the first term, from \eqref{hadamard} we have
\begin{equation}
\log\left|\frac{F(s)}{F(s+1)}\right|
=\sum_{\rho}\log\left|\frac{1-\frac{s}{\rho}}{1-\frac{s+1}{\rho}}\right|
=-\sum_{\rho}\log\left|\frac{s+1-\rho}{s-\rho}\right|.
\end{equation}
Thus,
\begin{equation}
\int_{1/2}^{3/2}\log\left|\frac{F(\sigma+it)}
{F(\sigma+it+1)}\right|d\sigma
=-\sum_{\rho}\int_{1/2}^{3/2}
\log\left|\frac{\sigma+it+1-\rho}{\sigma+it-\rho}\right|d\sigma.
\label{zerosum}
\end{equation}
Now, by the functional equation, the zeros of $\Lambda$ either lie on
the line $\Re(s)=\frac12$ or come in pairs $\rho, 1-\overline{\rho}$.
Applying Lemma \ref{hardlemma} with $w=\frac12+it-\rho$,
we see that \eqref{zerosum} is bounded below by
\begin{equation}
-(\log 4)\sum_{\rho}\Re\frac1{3/2+it-\rho}.
\end{equation}
Again by \eqref{hadamard}, this equals
\begin{equation}
\begin{aligned}
-(\log 4)\Re\frac{F'}{F}\!&\left(\frac32+it\right)\\
&=-(\log 4)\Re\left[
\frac{\gamma'}{\gamma}\!\left(\frac32+it\right)
+\frac{L'}{L}\!\left(\frac32+it\right)
+\frac{P'}{P}\!\left(\frac32+it\right)
+\frac{P'}{P}\!\left(\frac12+it\right)\right]\\
&\ge -(\log 4)\left[
\frac12\log\left|Q\!\left(\frac32+it\right)\right|
+\frac{4r}{\pi^2X^2}
+r\frac{z_{\theta}'}{z_{\theta}}\!\left(\frac32\right)
+\frac{2m}{X^2}\right].
\end{aligned}
\end{equation}

Altogether, we get
\begin{equation}
\begin{aligned}
\int_{1/2}^{\infty}\log|L(\sigma+it)|\,d\sigma
&\ge \left(\frac12-\log2\right)
\log\left|Q\!\left(\frac32+it\right)\right|\\
&+ r\left[\int_{3/2}^{5/2}\log z_{\theta}(\sigma)\,d\sigma
+\int_{3/2}^{\infty}\log z_{\theta}(\sigma)\,d\sigma
-(\log 4)\frac{z_{\theta}'}{z_{\theta}}\!\left(\frac32\right)\right]\\
&- r\left(\frac1{2\sqrt2X}+\frac{\frac4{\pi^2}\log(4e)+\frac14}{X^2}\right)
-\frac{4m\log{2}}{X^2}.
\end{aligned}
\end{equation}

Finally, we combine the upper bound for $t=t_2$ and lower bound for
$t=t_1$.  We get the stated main term plus error
\begin{equation}
r\left(\frac1{\sqrt2X}+\frac{\frac2{\pi^2}\log(16e^3)+\frac14}{X^2}\right)
+m\frac{4(\log2)+5/2}{X^2} < \frac{r}{\sqrt2(X-5)}.
\end{equation}
\end{proof}

\section{Rigorous Computation of $L$-functions}
\label{sec:rigorous}
The methods of Section \ref{turing} depend on a fast, rigorous algorithm
for evaluating $\Lambda(s)$.  We describe one such algorithm, based on
the Fast Fourier Transform, in this section.  We note that in the case
of the Riemann zeta function, a similar technique was developed and used
by Odlyzko and Sch\"onhage \cite{odlyzko-schonhage}.

Some algorithms for computing general $L$-functions were described by
Dokchitser \cite{dokchitser} and Rubinstein \cite{rubinstein}.
They ultimately boil down to the Cauchy integral formula:
\begin{equation}
\Lambda(s_0) = \frac1{2\pi i}\int\frac{\gamma(s)L(s)}{s-s_0}\,ds,
\end{equation}
where the contour consists of two vertical lines
enclosing $s_0$.  Writing $L(s)$ as a Dirichlet series and using the
functional equation, one is lead to study integrals of the form
\begin{equation}
\frac1{2\pi i}\int\frac{\gamma(s)n^{-s}}{s-s_0}\,ds,
\label{twoparam}
\end{equation}
taken along a vertical line far to the right.  Rubinstein, following an
idea of Lagarias and Odlyzko \cite{lagarias-odlyzko}, inserts a factor
designed to cancel the decay of the $\gamma$ factor,
e.g.\ $e^{i\frac{\pi r}4\eta s}$ for some $\eta$ close to $\pm 1$.
Without this factor, very high precision is required to calculate
$\Lambda(s)$ when $\Im(s)$ is large.

These algorithms are good when one is interested in computing
$\Lambda(s)$ at specific points, e.g.\ for locating zeros of $L(s)$
precisely.  They suffer from the disadvantage of being difficult to
carry out rigorously, basically because \eqref{twoparam} is a
two parameter family (indexed by $s_0$ and $n$) of integrals,
for which uniform asymptotics are hard to obtain in certain transition
ranges.

For Turing's method, we need an algorithm for {\em rigorously}
computing $\Lambda(s)$ for many values of $s$, not necessarily with
high precision.  For that we consider instead the one parameter
integrals
\begin{equation}
\frac1{2\pi i}\int\Lambda(s)e^{-zs}\,ds
\quad\mbox{and}\quad
\frac1{2\pi i}\int\gamma(s)e^{-zs}\,ds.
\end{equation}
These are essentially Fourier transforms, and they contain enough
information for evaluating $\Lambda(s)$ quickly, if one is
interested in many points.  They also involve only a single Mellin
transform, making rigorous computation more accessible.

Precisely, let $\eta\in(-1,1)$ and set
$F(t) := \Lambda\!\left(\frac12+it\right)
e^{\frac{\pi r}4\eta t}$.  Then the (inverse) Fourier transform of $F$ is
\begin{equation}
\begin{aligned}
\hat{F}(x)&:=\frac1{2\pi}\int_{-\infty}^{\infty}
F(t)e^{-ixt}\,dt
=\frac1{2\pi i}\int_{\Re(s)=\frac12}
\Lambda(s)e^{(x+i\frac{\pi r}4\eta)(1/2-s)}\,ds\\
&=\frac1{2\pi i}\int_{\Re(s)=2}
\Lambda(s)e^{(x+i\frac{\pi r}4\eta)(1/2-s)}\,ds
-\sum_{\rho\in\{1+\lambda_k:1\le k\le m\}}
\Res_{s=\rho}\Lambda(s)e^{(x+i\frac{\pi r}4\eta)(1/2-s)}.
\label{residuesum}
\end{aligned}
\end{equation}
The residue sum is straightforward to evaluate assuming we have complete
information on any poles of $L(s)$.
We multiply the Euler product \eqref{eulerprod} out
to a Dirichlet series, writing $L(s)=\sum_{n=1}^{\infty}a_nn^{-s}$.
Then the first term of \eqref{residuesum} is
\begin{equation}
\sum_{n=1}^{\infty}a_n\frac1{2\pi i}\int_{\Re(s)=2}
\gamma(s)e^{(x+i\frac{\pi r}4\eta)(1/2-s)}n^{-s}\,ds
=\epsilon\sum_{n=1}^{\infty}\frac{a_n}{\sqrt n}
G\!\left(x+\log\frac{n}{\sqrt N};\eta,\{\mu_j\}\right),
\label{Gseries}
\end{equation}
where
\begin{equation}
G(u;\eta,\{\mu_j\}) := \frac1{2\pi i}\int_{\Re(s)=2}
e^{(u+i\frac{\pi r}4\eta)(1/2-s)}
\prod_{j=1}^r\Gamma_{\R}(s+\mu_j)
\,ds.
\label{Gdef}
\end{equation}

Let us assume for now that we have a procedure to compute
$G(u;\eta,\{\mu_j\})$,
and thereby $\hat{F}(x)$, to prescribed precision; we return to this
point in Section \ref{sec:Gcompute} below.
In order to use the FFT to compute $F$ from $\hat{F}$, we first need to
discretize the problem.  To that end,
let $A,B>0$ be parameters such that $q=AB$ is an integer.  By the Poisson
summation formula,
\begin{equation}
\sum_{k\in\Z}F\!\left(\frac{m}A+kB\right)
=\frac{2\pi}{B}\sum_{k\in\Z}
\hat{F}\!\left(\frac{2\pi k}B\right)\!e\!\left(\frac{km}{AB}\right)
=\frac{2\pi}{B}
\sum_{n\,(\mbox{\scriptsize mod }q)}
e\!\left(\frac{mn}q\right)\!
\sum_{k\in\Z}\hat{F}\!\left(\frac{2\pi n}B+2\pi Ak\right)\!.
\end{equation}
Thus, the functions
$\widetilde{F}(m):=\sum_{k\in\Z}F\!\left(\frac{m}A+kB\right)$
and
$\widetilde{\hat{F}}(n):=
\sum_{k\in\Z}\hat{F}\!\left(\frac{2\pi n}B+2\pi Ak\right)$,
which are periodic in $m,n$ with period $q$, form a discrete Fourier
transform pair.

Note that since $F$ is real-valued, $\hat{F}(-x) = \overline{\hat{F}(x)}$.
Thus, for $|n|\le q/2$ we have
\begin{equation}
\widetilde{\hat{F}}(n)=\hat{F}\!\left(\frac{2\pi n}B\right)
+\sum_{k=1}^{\infty}
\hat{F}\!\left(\frac{2\pi n}B+2\pi kA\right)
+\sum_{k=1}^{\infty}
\overline{\hat{F}\!\left(-\frac{2\pi n}B+2\pi kA\right)}.
\label{hatFsmall}
\end{equation}
For $A$ even moderately large, the terms for $k\ge 1$ fall within the
asymptotic range.  Precise bounds are given in Section
\ref{sec:asymp} below; in particular, we may apply Lemma \ref{lem:Fasymp}
with $x=2\pi(A\pm\frac{n}{B})$ to compute the sums over $k$.  Hence, it
suffices to calculate $\hat{F}(2\pi n/B)$ for $0\le n\le q/2$.
On the other hand, to compute $F(m/A)$, we need to bound the terms of
$\widetilde{F}(m)$ for $k\ne 0$.  We have already obtained a
suitable bound for the $L$-function in Lemma \ref{Lupperbound}.
The sum of this bound over $k\ne 0$ is the content of Lemma
\ref{lem:Lambdasum}.

\subsection{Computing $G(u;\eta,\{\mu_j\})$} \label{sec:Gcompute}
For brevity, some of the results of this section are only sketched.
Our emphasis is on the details necessary for
rigorous computation.  For more general background
information we refer the reader to
\cite{booker2,dokchitser,rubinstein}.

One simple method for calculating integrals such as \eqref{Gdef}
that is easy to make rigorous is the power and log
series, obtained by shifting the contour of \eqref{Gdef} to the left:
\begin{equation}
\begin{aligned}
G(u;\eta,\{\mu_j\}) &= \sum_{\rho\in\C}
\Res_{s=\rho}\left(
e^{(u+i\frac{\pi r}4\eta)(1/2-s)}
\prod_{j=1}^r\Gamma_{\R}(s+\mu_j)
\right)\\
&= \sum_{\mathrm{poles}\,\rho}
P(u;\rho,\eta,\{\mu_j\})e^{(1/2-\rho)u},
\end{aligned}
\label{powerlog}
\end{equation}
where $P(u;\rho,\eta,\{\mu_j\})$ is a polynomial of degree one less
than the order of the pole at $\rho$.

For example, in the case of Galois representations,
the $\mu_j$ are all either $0$
or $1$, and the residues in \eqref{powerlog} may be evaluated by the
following:
\begin{equation}
\begin{aligned}
&\begin{aligned}
\Gamma_{\R}(s)=
\frac2{s+2k}\frac{(-2\pi)^k}{(2k)!!}
\exp&\!\left[
\!\left(\sum_{n=1}^k\frac1{2n}-
\frac12\log(\pi e^{\gamma})\right)\!(s+2k)\right.\\
&+\left.\sum_{j=2}^{\infty}\frac{1}j
\!\left((-1)^j2^{-j}\zeta(j)+\sum_{n=1}^k\frac1{(2n)^j}\right)\!
(s+2k)^j\right]
\end{aligned}\\
&\hspace{11mm}\begin{aligned}
=\frac{(-2\pi)^k}{(2k-1)!!}
\exp&\!\left[
\!\left(\sum_{n=1}^k\frac1{2n-1}-
\frac12\log(4\pi e^{\gamma})\right)\!(s+2k-1)\right.\\
&+\left.\sum_{j=2}^{\infty}\frac{1}j
\!\left((-1)^j(1-2^{-j})\zeta(j)+\sum_{n=1}^k\frac1{(2n-1)^j}\right)\!
(s+2k-1)^j\right],
\end{aligned}
\end{aligned}
\end{equation}
for any integer $k\ge 0$.  For general $\mu_j$, we need an algorithm to
calculate the values of $\Gamma$ and its derivatives at an arbitrary
point in the complex plane; we assume without further comment that this
is available when necessary.

We say that $\mu_j$ and $\mu_k$ are equivalent if $\mu_j-\mu_k\in 2\Z$.
For $\mu$ ranging over an equivalence class, the functions
$\Gamma_{\R}(s+\mu)$ share all but finitely many poles.  Thus,
\eqref{powerlog} may be broken naturally into parts corresponding to each
class.  We can bound the tail of each part as follows.
\begin{lemma}
Let $\rho$ be a pole of 
$g(s)=e^{(u+i\frac{\pi r}4\eta)(1/2-s)}
\prod_{j=1}^r\Gamma_{\R}(s+\mu_j)$ of order $n$,
with $\Re(\rho+\mu_j)\le 0$ for $j=1,\ldots,r$ and
$(2\pi)^re^{2u}<\frac12\prod_{j=1}^r(|2-\rho-\mu_j|-1)$.
Let $c_j$ be the coefficients of the polar part of $g$
around $\rho$, i.e.\ such that $g(s+\rho)-\sum_{j=1}^nc_js^{-j}$ is
holomorphic at $s=0$.  Then
\begin{equation}
\left|\sum_{k=1}^{\infty}\Res_{s=\rho-2k}g(s)\right|
< \max|c_j|.
\end{equation}
\end{lemma}
\begin{proof}
First note that
\begin{equation}
\begin{aligned}
g(s+\rho-2)&=(2\pi)^re^{2u+i\frac{\pi r}2\eta}g(s+\rho)
\prod_{j=1}^r(s+\rho-2+\mu_j)^{-1}\\
&=\frac{(-2\pi)^re^{2u+i\frac{\pi r}2\eta}}
{\prod_{j=1}^r(2-\rho-\mu_j)}
\cdot\frac{g(s+\rho)}
{\prod_{j=1}^r\bigl(1-\frac{s}{2-\rho-\mu_j}\bigr)}.
\end{aligned}
\label{eq:gminus2}
\end{equation}

Next, let $f(s)$ be a meromorphic function with polar part
$a_1s^{-1}+\ldots+a_ns^{-n}$ at $0$.  If $x$ is a complex number with
$|x|<1$ then the function $\frac{f(s)}{1-xs}$ has polar part
$a_1's^{-1}+\ldots+a_n's^{-n}$, where
$a_j'=\sum_{k=0}^{n-j}a_{j+k}x^k$.  Thus,
\begin{equation}
\max|a_j'| \le \frac{\max|a_j|}{1-|x|}.
\label{eq:xsmax}
\end{equation}

Let $c_1's^{-1}+\ldots+c_n's^{-n}$ be the polar part of $g(s+\rho-2)$.
Applying \eqref{eq:xsmax} $r$ times, we see from \eqref{eq:gminus2} that
\begin{equation}
\begin{aligned}
\max|c_j'|& \le
\frac{(2\pi)^re^{2u}}{\prod_{j=1}^r|2-\rho-\mu_j|}
\cdot\frac{\max|c_j|}{\prod_{j=1}^r\bigl(1-\frac1{|2-\rho-\mu_j|}\bigr)}\\
&=\frac{(2\pi)^re^{2u}\max|c_j|}{\prod_{j=1}^r(|2-\rho-\mu_j|-1)}
< \frac12\max|c_j|.
\end{aligned}
\end{equation}
Repeating this procedure, we see that the coefficients of the polar part
of $g(s+\rho-2k)$ are $< 2^{-k}\max|c_j|$.  The conclusion follows.
\end{proof}

The lemma says roughly that if we compute the residue sum for all poles
with real part down to $\Re(\rho)$, the tail of the series (from poles
at $\rho-2k$) may be bounded by the data from the last term added.
Moreover, \eqref{eq:gminus2} gives an algorithm for computing the data
at $\rho-2$ from that at $\rho$, and shows that the terms eventually
decrease factorially.  Thus, we may use this to compute
$G(u;\eta,\{\mu_j\})$ to any desired precision for a given $u$.

Since the $\mu_j$ are arbitrary, this procedure is general enough to
compute derivatives of $G(u;\eta,\{\mu_j\})$ as well.  For instance,
for any $k$ we have
\begin{equation}
G'(u;\eta,\{\mu_j\})
=(\mu_k+1/2)G(u;\eta,\{\mu_j\})-2\pi G(u;\eta,\{\mu_j'\}),
\end{equation}
where $\mu_j' = \mu_j$ if $j\ne k$ and $\mu_k' = \mu_k+2$.
Higher derivatives may be computed in a similar fashion.  (In fact $G$
satisfies an $r$th order differential equation, due to the recurrence
for $\Gamma$; thus, the derivatives of all orders are determined from
the first $r$.)

Note that for $u$ large this method requires high precision due to
cancellation, and is therefore inefficient.  The essential point that
makes it worthwhile is that for a given $\gamma$ factor the
calculations need only be performed once, as one can develop local
approximations to $G(u;\eta,\{\mu_j\})$ for later rapid evaluation.
The computation may then be recycled and used for any $L$-function with
the same $\mu_j$; this is useful for functions in an arithmetic family,
such as Artin $L$-functions.

More precisely, suppose we wish to calculate $G(u;\eta,\{\mu_j\})$ for $u$
in an interval $I$.  Choose $\varepsilon>0$ and sample points $u_m$ such
that each $u\in I$ is contained in a unique interval
$[u_m-\varepsilon,u_m+\varepsilon)$.  For $u$ in the $m$th interval,
we have by Taylor's theorem
\begin{equation}
G(u;\eta,\{\mu_j\}) = \sum_{k=0}^{K-1}\frac{G^{(k)}(u_m;\eta,\{\mu_j\})}{k!}
(u-u_m)^k +
\Theta\!\left(\max_{|u^*-u_m|\le\varepsilon}
\frac{\bigl|G^{(K)}(u^*;\eta,\{\mu_j\})\bigr|}{K!}\varepsilon^K\right).
\label{localapprox}
\end{equation}
We may evaluate the derivatives precisely using \eqref{powerlog}.  As
for the $K$th derivative, a uniform bound is obtained by shifting the
contour of \eqref{Gdef} to $\Re(s)=\frac12$:
\begin{equation}
\frac{\bigl|G^{(K)}(u;\eta,\{\mu_j\})\bigr|}{K!} \le
\frac1{2\pi}\int_{-\infty}^{\infty}\frac{|t|^K}{K!}
e^{\frac{\pi r}4\eta t}
\prod_{j=1}^r\left|\Gamma_{\R}\!\left(\frac12+it+\mu_j\right)\right|
\,dt.
\end{equation}
For large $K$ this is of size $\left(\frac4{\pi r(1-\eta)}\right)^K$;
thus as long as $\varepsilon$ is small compared to $\frac{\pi
r}4(1-\eta)$, we may compute and store the coefficients of
\eqref{localapprox}, yielding a fast method to calculate
$G(u;\eta,\{\mu_j\})$ for any $u\in I$.

Moreover, we can improve the efficiency of our algorithm if the sample
points of \eqref{hatFsmall} coincide with multiplies of $2\varepsilon$,
i.e.\ if $\frac{\pi}{\varepsilon B}\in\Z$.  For any given sample point
$x$, we approximate ${\hat F}(x)$ via a truncated series (keeping track
of the error terms from \eqref{localapprox} and Lemma \ref{lem:Ftail}
below):
\begin{equation}
\begin{aligned}
&\hspace{-1cm}\sum_{n=1}^M\frac{a_n}{\sqrt n}
G\!\left(x+\log\frac{n}{\sqrt N};\eta,\{\mu_j\}\right)\\
&\approx
\sum_m\sum_{\log\frac{n}{\sqrt N}\in [u_m-\varepsilon,u_m+\varepsilon)}
\frac{a_n}{\sqrt n}\sum_{k=0}^{K-1}
\frac{G^{(k)}(x+u_m;\eta,\{\mu_j\})}{k!}
\left(\log\frac{n}{\sqrt N}-u_m\right)^k\\
&=
\sum_{k=0}^{K-1}\sum_m
\frac{G^{(k)}(x+u_m;\eta,\{\mu_j\})}{k!}
S^{(k)}_m,
\end{aligned}
\label{eq:Fxlocal}
\end{equation}
where
\begin{equation}
S^{(k)}_m :=
\sum_{\log\frac{n}{\sqrt N}\in [u_m-\varepsilon,u_m+\varepsilon)}
\frac{a_n}{\sqrt n}
\left(\log\frac{n}{\sqrt N}-u_m\right)^k.
\end{equation}
Since $x+u_m$ is another sample point $u_{m'}$, the $k$th term of
\eqref{eq:Fxlocal} is a convolution of the sequences (indexed by $m$)
$\frac{G^{(k)}(u_m;\eta,\{\mu_j\})}{k!}$ and $S^{(k)}_m$;
thus, we may evaluate it efficiently for all $x$ simultaneously by
appealing again to the FFT.

\subsection{Complexity}
We may now consider the complexity of the algorithm.  Note that by
Stirling's formula, $F(t)$ decays roughly like $e^{-(1-\eta)\frac{\pi
r}4t}$ for $t>0$.  Ideally we should choose $1-\eta$ of size $T^{-1}$
in order to compute values up to height $T$.  Adjusting the constant of
proportionality (i.e.\ choosing $\eta$ relatively close to or far from
$1$) allows us to trade off the computational precision and number of
coefficients needed to overcome the error terms below.  Finding a good
compromise between these two is best done by trial and error; cf.\ Section
\ref{numerics}.

Since $\delta\asymp T^{-1}$, Lemma \ref{lem:Ftail} shows that in order
to compute \eqref{eq:Fxlocal} we need on the order of $\sqrt{NT^r}$
terms of \eqref{Gseries}, or roughly the square root of the analytic
conductor.  Note that the values of $A$ and $B$ enter only in the
Fourier transforms, and do not significantly affect the computation of
\eqref{eq:Fxlocal}.  We set $B$ equal to a multiple of $T$, depending
on the chosen value of $\eta$.  As for $A$, as mentioned in Section
\ref{sec:explicit}, the density of zeros of $F(t)$ around height $T$
is $\frac1{2\pi}\log N\!\left(\frac{T}{2\pi}\right)^r$; one can expect to
take $A$ equal to a multiple of this.  Thus, this method has complexity
consistent with computing a {\em single} value by the approximate
functional equation, after which we get many values in mean time
$O_{\varepsilon}\bigl((NT)^{\varepsilon}\bigr)$, which is essentially
best possible.  The gain comes from the fact, as emphasized above,
that only a single $G$-function is involved.

\subsection{Asymptotics} \label{sec:asymp}
To complete our understanding of $G(u;\eta,\{\mu_j\})$, in order to
accurately calculate \eqref{Gseries}, we need an asymptotic bound for
large $u$.  If we write $\mu=-\frac12+\frac1r\bigl(
1+\sum_{j=1}^r\mu_j\bigr)$ then by the method of stationary phase, we have
\begin{equation}
G(u;\eta,\{\mu_j\})=
\sqrt{\frac{2^{r+1}}r}e^{\mu(u+i\frac{\pi r}4\eta)}
\exp\!\left(-e^{\frac2r(u+i\frac{\pi r}4\eta)}\right)\!
\bigl(1+O\bigl(e^{-2u/r}\bigr)\bigr),
\label{Gasymp}
\end{equation}
where the implied constant depends on the $\mu_j$.
For $r=1$ the formula is exact, i.e.\ the $O$ term is $0$.  For $r>1$,
one can work out explicit constants case by case, which
is preferable if sharp error terms are desired.  Otherwise, we get a
bound that is close to \eqref{Gasymp} simply by shifting the contour of
\eqref{Gdef} to the right.
\begin{lemma}
\label{lem:Gasymp}
Let $\delta=\frac{\pi}2(1-|\eta|)$,
$\nu_j=\frac{\Re(\mu_j)-1}2+\frac1{2r}$,
$\mu=-\frac12+\frac1r\bigl(1+\sum_{j=1}^r\mu_j\bigr)$,
$K=2\sqrt{\frac{2^{r+1}}r\frac{e^{\delta(r-1)}}{\delta}}
e^{-\frac{\pi r\eta\Im(\mu)}4}$,
and $X=\pi r\delta e^{-\delta}e^{2u/r}$.
Then for $X\ge r$,
\begin{equation}
\bigl|G(u;\eta,\{\mu_j\})\bigr| \le
Ke^{\Re(\mu)u}e^{-X}
\prod_{j=1}^r\left(1+\frac{r\nu_j}X\right)^{\nu_j}.
\label{Gbound}
\end{equation}
\end{lemma}
\begin{remark}
This is within a factor $O\bigl(\delta^{-1/2}\bigr)$ of the correct
asymptotic if $\delta\ll re^{-u/r}$.
\end{remark}
\begin{proof}
We write $s=2\sigma+2it$ in \eqref{Gdef} to get
\begin{equation}
|G(u;\eta,\{\mu_j\})|\le
\pi^{-r(\sigma+\frac{\Re(\mu)}2+\frac14)-\frac12}
e^{u(\frac12-2\sigma)-\frac{\pi r\eta\Im(\mu)}4}
\int_{-\infty}^{\infty}\prod_{j=1}^r
\left|\Gamma\!\left(\sigma+it+\frac{\mu_j}2\right)
e^{\frac{\pi \eta}4(2t+\Im(\mu_j))}
\right|dt.
\label{Gbound1}
\end{equation}

Applying H\"older's inequality, we get integrals of the form
\begin{equation}
\int_{-\infty}^{\infty}|\Gamma(a+it)|^re^{\frac{\pi r\eta}2 t}\,dt,
\end{equation}
where $a=\sigma+\Re(\mu_j)/2$.  Assuming $\sigma\ge 1$, we may apply the
inequality
\begin{equation}
|\Gamma(a+it)|\le \sqrt{2\pi}(a+|t|)^{a-1/2}e^{-\pi|t|/2}
\quad\mbox{for }a\ge\frac12.
\end{equation}
(To see this, note that
$\bigl|\Gamma\bigl(\frac12+it\bigr)\bigr|=\sqrt{\pi {\rm sech}\,\pi t}$,
use the recurrence for $\Gamma$ and \eqref{Gammalogderiv} to reduce to
the region $\frac12\le a\le \frac32$, $0\le t\le 2$, where the inequality
may be checked computationally.)  Thus, we have
\begin{equation}
\begin{aligned}
\int_{-\infty}^{\infty}|\Gamma(a+it)|^re^{\frac{\pi r\eta}2 t}\,dt
&\le 2(2\pi)^{r/2}\int_0^{\infty}(a+t)^{r(a-1/2)}e^{-r\delta t}\,dt\\
&\le 2(2\pi)^{r/2}e^{\delta ra}
\frac{\Gamma(r(a-1/2)+1)}{(\delta r)^{r(a-1/2)+1}} \\
&\le 2(2\pi)^{\frac{r+1}2}\sqrt{\frac{e^{\delta(r-1)}}{\delta r}}
\left(\frac{e^{\delta}}{\delta}\frac{a-1/2+1/2r}e\right)^{r(a-1/2)+1/2}.
\end{aligned}
\end{equation}
Substituting this bound into \eqref{Gbound1} and collecting terms we obtain
\begin{equation}
Ke^{\Re(\mu)u}\prod_{j=1}^r
\left(\frac{\sigma+\nu_j}{eX/r}\right)^{\sigma+\nu_j}
\le
Ke^{\Re(\mu)u}e^{-\sigma r}\prod_{j=1}^r
\left(\frac{r\sigma}{X}\right)^{\sigma+\nu_j}
\left(1+\frac{\nu_j}{\sigma}\right)^{\nu_j}.
\end{equation}
The result follows upon taking $\sigma=X/r$.
\end{proof}

With this bound in hand, we can estimate the error in truncating the
series \eqref{Gseries}.
\begin{lemma}
\label{lem:Ftail}
Let $M$ be a positive integer, $x\in\R$.  Let $\delta,\nu_j,\mu,K$
be as in Lemma \ref{lem:Gasymp} and set
$X=\pi r\delta e^{-\delta}\bigl(e^x/\sqrt N\bigr)^{2/r}$.
Let $C,\alpha\ge0$ be such that $|a_n|\le Cn^{\alpha}$ for all $n$,
and put $c=\Re(\mu)+\frac12+\alpha$, $c'=\max(cr/2-1,0)$.
Then for $XM^{2/r}>\max(c',r)$,
\begin{equation}
\begin{aligned}
&\hspace{-1cm}\left|\sum_{n>M}\frac{a_n}{\sqrt n}
G\!\left(x+\log\frac{n}{\sqrt N};\eta,\{\mu_j\}\right)\right|\\
&\le \frac{Kr}2\left(\frac{e^x}{\sqrt N}\right)^{\Re(\mu)}
\frac{CM^ce^{-XM^{2/r}}}{XM^{2/r}-c'}
\prod_{j=1}^r
\left(1+\frac{r\nu_j}{XM^{2/r}}\right)^{\nu_j}.
\end{aligned}
\end{equation}
\end{lemma}
\begin{remark}
Different values of $C$ and $\alpha$ are appropriate for different
ranges.  For small $M$, one can take $C=1,\alpha=\log_2r+\theta$,
while for larger $M$ it is better to choose a smaller value of $\alpha$
and compute $C$ from the coefficients.
\end{remark}
\begin{proof}
Using Lemma \ref{lem:Gasymp}, we have
\begin{equation}
\sum_{n>M}\frac{|a_n|}{\sqrt n}
\left|G\!\left(x+\log\frac{n}{\sqrt N};\eta,\{\mu_j\}\right)\right|
\le K'\sum_{n>M}n^{c-1}e^{-Xn^{2/r}},
\label{eq:Ftail1}
\end{equation}
where
\begin{equation}
K'=CK\left(\frac{e^x}{\sqrt N}\right)^{\Re(\mu)}
\prod_{j=1}^r\left(1+\frac{r\nu_j}{XM^{2/r}}\right)^{\nu_j}.
\end{equation}
The condition on $X$ ensures that the terms of \eqref{eq:Ftail1}
are monotonically decreasing.  Thus, we can estimate by the integral
\begin{equation}
K'\int_M^{\infty}t^{c-1}e^{-Xt^{2/r}}\,dt
=\frac{K'r}2\bigl(X^{-r/2}\bigr)^c\int_{XM^{2/r}}^{\infty}
y^{cr/2-1}e^{-y}\,dy
\le\frac{K'r}2\frac{M^ce^{-XM^{2/r}}}{XM^{2/r}-c'}.
\end{equation}
\end{proof}

The next two lemmas bound the error introduced in discretization.
\begin{lemma}
\label{lem:Fasymp}
Let $x\in\R$, $A\ge \frac1{2\pi}$, and let notation be as in Lemma
\ref{lem:Ftail}.  Then for $X>\max(c',r)$,
\begin{equation}
\begin{aligned}
\sum_{k=0}^{\infty}\hat{F}(x+2\pi kA)&=
-\sum_{\rho\in\{1+\lambda_k:1\le k\le m\}}
\Res_{s=\rho}\frac{\Lambda(s)e^{(x+i\frac{\pi r}4\eta)(1/2-s)}}
{1-e^{2\pi A(1/2-s)}}\\
&+\Theta\!\left[
\frac{K}{1-e^{-\pi A}}
\!\left(\frac{e^x}{\sqrt N}\right)^{\Re(\mu)}e^{-X}
\!\left(1+\frac{Cr/2}{X-c'}\right)
\prod_{j=1}^r\left(1+\frac{r\nu_j}X\right)^{\nu_j}\right]\!.
\end{aligned}
\end{equation}
\end{lemma}
\begin{proof}
The residue sum comes from summing the polar part of \eqref{residuesum}
with $x+kA$ in place of $x$.  For the rest, we apply Lemmas
\ref{lem:Gasymp} and \ref{lem:Ftail} (with $M=1$) to get the bound
\begin{equation}
K\!\left(\frac{e^x}{\sqrt N}\right)^{\Re(\mu)}\!e^{-X}
\!\left(1+\frac{Cr/2}{X-c'}\right)
\prod_{j=1}^r
\!\left(1+\frac{r\nu_j}X\right)^{\nu_j}
\end{equation}
for the $k=0$ term.  To pass from this to the $k$th term, we multiply by
a factor not exceeding
\begin{equation}
\begin{aligned}
e^{2\pi kA\Re(\mu)}e^{-X(\exp(4\pi kA/r)-1)}
&=\exp\!\left(2\pi kA\left[\Re(\mu)-\frac{2X}r-
\frac{X}{2\pi kA}\left(e^{4\pi kA/r}-1-\frac{4\pi kA}r\right)\right]\right)\\
&\le \exp\!\left(-2\pi kA\left[\alpha+\frac12-\frac2r
+\frac{4\pi kA}r\right]\right)
\le e^{-\pi kA}.
\end{aligned}
\end{equation}
The result follows on summing the geometric series.
\end{proof}

\begin{lemma}
\label{lem:Lambdasum}
Let $t\in\R$ and put $s=\frac12+it$,
\begin{equation}
E=Z_{\theta}(3/2)^r
|\gamma(s)|e^{\frac{\pi r}4\eta t}
\left|Q(s)\frac{P(s+1)^2P(s-2)}{P(s)^2P(s-1)}\right|^{1/2},
\end{equation}
and
\begin{equation}
\beta = \frac{\pi r}4
-\frac12\sum_{j=1}^r\arctan\frac{\Re(s+\mu_j)}{|\Im(s+\mu_j)|}
-\frac4{\pi^2}\sum_{j=1}^r
\frac1{|\Im(s+\mu_j)^2-\Re(s+\mu_j)^2|}.
\end{equation}
\begin{enumerate}
\item
If $\Im(s+\mu_j) > 0$ for all $j=1,\ldots,r$ and
$\beta-\frac{\pi r}4\eta > 0$ then
\begin{equation}
\left|\sum_{k=0}^{\infty}F\bigl(t+kB\bigr)\right|
\le \frac{E}{1-e^{-(\beta-\frac{\pi r}4\eta)B}}.
\end{equation}
\item
If $\Im(s+\mu_j) < 0$ for all $j=1,\ldots,r$ and
$\beta+\frac{\pi r}4\eta > 0$ then
\begin{equation}
\left|\sum_{k=0}^{\infty}F\bigl(t-kB\bigr)\right|
\le \frac{E}{1-e^{-(\beta+\frac{\pi r}4\eta)B}}.
\end{equation}
\end{enumerate}
\end{lemma}
\begin{proof}
We treat only the first case, the second being similar.  Lemmas
\ref{Lupperbound} and \ref{zetalemma} imply the bound $|F(t)|\le E$.
We consider the same bound with $t$ replaced by $t+kB$.  Note that if
$|\Im(s+\mu_j)|$ increases for all $j$ then the factor involving $P$ is
non-increasing.  For the $\gamma$ and $Q$ factors, by the mean value
theorem we have
\begin{equation}
\log\!\left(\left|\frac{\gamma(s+ikB)}{\gamma(s)}\right|
\left|\frac{Q(s+ikB)}{Q(s)}\right|^{1/2}\right)
=-kB\Im\!\left(
\frac{\gamma'}{\gamma}(s^*)+\frac12\frac{Q'}{Q}(s^*)\right)
\end{equation}
for some $s^*$ on the line between $s$ and $s+ikB$.
Using \eqref{Gammalogderiv}, this is
\begin{equation}
\begin{aligned}
-kB&\Im
\sum_{j=1}^r\left(
\frac12\log\frac{s^*+\mu_j}2
+\Theta\!\left(\frac{4/\pi^2}{|\Im(s^*+\mu_j)^2-\Re(s^*+\mu_j)^2|}
\right)\right)\\
&\le
-kB\!\left(
\frac{\pi r}4
-\frac12\sum_{j=1}^r\arctan\frac{\Re(s^*+\mu_j)}{\Im(s^*+\mu_j)}
-\frac4{\pi^2}\sum_{j=1}^r
\frac1{|\Im(s^*+\mu_j)^2-\Re(s^*+\mu_j)^2|}\right)\\
&\le -\beta kB.
\end{aligned}
\end{equation}
Thus, $|F(t+kB)|\le Ee^{-(\beta-\frac{\pi r}4\eta)kB}$.
The conclusion follows.
\end{proof}

\section{Numerical results} \label{numerics}
We have applied the methods described in the previous sections to a
few examples of splitting fields of polynomials with Galois group $S_5$
and $A_5$, as listed in Table \ref{tab:examples}.  For the $A_5$ cases,
the Artin conjecture is true for all representations by known cases
of functoriality \cite{kiming,jehanne2}.  That speeds up the process,
since we may apply Turing's method to the Artin $L$-functions directly.
For the $S_5$ examples we verify both conjectures.  As expected, we found
no counterexamples to either conjecture in the tested range $|t|\le 100$.

\begin{table}[h]
\begin{tabular}{|r|r|l|}
\hline
polynomial & group & splitting field discriminant \\ \hline
$x^5-68x-68$ & $S_5$ & $2^{96}3^{60}17^{96}$ \\
$x^5-x^4-8x^3+10x^2-x-5$ & $S_5$ & $2^{160}3^{96}7^{96}$ \\
$x^5-x^4+3x^3-11x^2-8x-8$ & $S_5$ & $2^{220}13^{96}$ \\
$x^5+2x^3-4x^2-2x+4$ & $A_5$ & $2^{90}73^{30}$ \\
$x^5+20x+16$ & $A_5$ & $2^{90}5^{78}$ \\
$x^5-x^4+8x^3-6x^2+14x-6$ & $A_5$ & $2^{90}193^{30}$ \\
$x^5-7x^3-17x^2+18x+73$ & $A_5$ & $2^{40}487^{30}$ \\
$x^5+8x^3+7x^2+172x+53$ & $A_5$ & $2083^{30}$ \\ \hline
\end{tabular}
\caption{Tested polynomials}
\label{tab:examples}
\end{table}

To illustrate the methods, we discuss in detail
the $S_5$ field of discriminant $2^{96}3^{60}17^{96}$ given by
the polynomial $f(x)=x^5-68(x+1)$.  Recall that $S_5$ has seven
irreducible representations.  We label them $1$, $\chi$, $\rho_4$,
$\rho_4'=\rho_4\otimes\chi$, $\rho_5$, $\rho_5'=\rho_5\otimes\chi$
and $\rho_6$, where $\chi$ is the sign character and the subscripts
indicate the dimensions.  As it will turn out, the limiting factor in
our computations is the conductor of $\rho_6$, which in our example
is $36081072=2^4 3^3 17^4$.  This is the smallest among the table of
$S_5$ polynomials given in \cite{nfdatabase}; since that
table is ordered by the conductor of $\rho_4$, it is likely that smaller
examples exist.  (We note, however, that if one is interested only in
verifying some instances of Artin's conjecture and not the Riemann
hypothesis, the holomorphy of $L(s,\rho_5')$ may be checked much more
easily; there the limiting factor is the conductor of $\rho_4$, of which
\cite{nfdatabase} yields examples as small as $1609$.  We have not
pursued this possibility.)

Note that $1$, $\chi$ and $\rho_6$ are monomial representations, so
Artin's conjecture is true for those.  Equation \eqref{eq:zetaratio}
below shows that $L(s,\rho_4)$ and $L(s,\rho_5)$ are holomorphic
except possibly at the zeros of $\zeta(s)$.  Twisting by $\chi$, we
see similarly that $L(s,\rho_4')$ and $L(s,\rho_5')$ are holomorphic
away from the zeros of $L(s,\chi)$.  Moreover, we learn from GAP
that the representations $\rho_4\oplus\rho_6$, $\rho_4'\oplus\rho_6$,
$\rho_5\oplus\rho_4'\oplus\chi$ and $\rho_5'\oplus\rho_4\oplus1$ are
all monomial.  Thus, in order to verify the holomorphy of $L(s,\rho_4)$,
$L(s,\rho_4')$, $L(s,\rho_5)$ and $L(s,\rho_5')$ it is enough to check
that $L(s,\rho_6)$ and $L(s,\rho_4'\oplus\chi)$ are non-vanishing
at zeros of $\zeta(s)$ and, similarly, that $L(s,\rho_6)$ and
$L(s,\rho_4\oplus1)$ do not vanish at the zeros of $L(s,\chi)$.  Applying
Turing's method to these functions as well as $L(s,\rho_5\oplus1)$ and
$L(s,\rho_5'\oplus\chi)$, we can deduce both Artin's conjecture and the
Riemann hypothesis (up to the tested height) for all representations.

In what follows we describe the numerical procedure in detail for
$L(s,\rho_6)$.  First we must choose a value of $\eta$ to use for the
computation of $G(u;\eta,\{\mu_j\})$.  One can aim to limit either the
number of Dirichlet coefficients $a_n$ or the precision required in the
computation.  Since the coefficients are relatively easy to compute in
our case (we have $2^{32}$ of them), we try for the latter.  The largest
error comes from Lemma \ref{lem:Ftail} with $x=0$, and is of size roughly
$M^c\exp\bigl(-\pi r\delta e^{-\delta}(M/\sqrt{N})^{2/r}\bigr)$, where
$\delta=\frac{\pi}2(1-|\eta|)$, $r=6$, $N=36081072$ and $M=2^{32}$.
Examining the local factors at small primes we determine that
$|a_n|\le 1.26n^{\log_{4243}6}$, yielding $c=\frac23+\log_{4243}6$.
(The $\mu_j$ in this case are $0,0,0,1,1,1$.)  This error term should be
compared to the size of the function being evaluated, which is roughly
$e^{-\delta rt/2}$.  From Theorem \ref{Stbound} we find that to apply
Turing's method up to height $t=100$ we need to be able to compute
the $L$-function up to about $t=115$.  Trying a few values of $\eta$,
we find that with $\eta=0.98$ the error terms are of size $10^{-14}$,
compared to $10^{-5}$ for the size of the function.
Thus, with this choice we should use a precision of at least $14$ digits;
in fact we carry out most computations to $30$ digits.

All computations were performed on a 3GHz PC running Linux.
They were divided into several steps:
\begin{enumerate}
\item Computing the Dirichlet coefficients $a_n$;
\item Estimating zeros by the explicit formula;
\item Computing $G^{(k)}(u_m;\eta,\{\mu_j\})$;
\item Computing $S^{(k)}_m$;
\item Computing \eqref{eq:Fxlocal} and $L(s,\rho_6)$ by FFT;
\item Turing's method.
\end{enumerate}
To ensure correct results, we used the arbitrary precision interval
arithmetic package MPFI \cite{MPFI} for steps 3 through 7.  We discuss
the steps in more detail below.

\subsection{Computing $a_n$} \label{numerics:coeff}
First we consider methods of computing the coefficients of the
$L$-functions $L(s,\rho)$ for all irreducible representations $\rho$.
One way is to express them as ratios of products of Hecke $L$-functions,
as given by Brauer's theorem; in fact for $S_5$ we may express each in
terms of Dedekind zeta functions of intermediate fields.  Precisely,
let $k$ be the quadratic extension of $\Q$ associated to $\chi$,
$F=\Q(x_1)\subset M=\Q(x_1,x_2)$ where $x_1$ and $x_2$ are distinct roots
of $f$, and $E=\Q(y)\subset K$ where $y$ is a root of the sextic resolvent
(a formula for which is given in \cite{dummit}); then we have
\begin{equation}
\begin{aligned}
L(s,\chi) &= \frac{\zeta_k(s)}{\zeta(s)},&\quad
L(s,\rho_4) &= \frac{\zeta_F(s)}{\zeta(s)},&\quad
L(s,\rho_5) &= \frac{\zeta_E(s)}{\zeta(s)},&\\
L(s,\rho_6)&=
\frac{\zeta_k(s)\zeta_E(s)\zeta_M(s)}{\zeta_{kE}(s)\zeta_F(s)^2},&\quad
L(s,\rho_4')&=
\frac{\zeta(s)\zeta_{kF}(s)}{\zeta_k(s)\zeta_F(s)},&\quad
L(s,\rho_5')&=
\frac{\zeta(s)\zeta_{kE}(s)}{\zeta_k(s)\zeta_E(s)}.&
\end{aligned}
\label{eq:zetaratio}
\end{equation}
In turn, we may compute each of the Dedekind zeta functions
using the ideal factorization functions built in to PARI \cite{PARI}.
This facilitates the computation of local factors at primes dividing
the discriminant, allowing us to avoid a detailed study of the possible
types of ramification.  However, it is not well-suited to working out
many coefficients.

Fortunately, there is a faster method that works for all but finitely
many primes.  Table \ref{tab:S5factors} shows the unramified local
factors for each representation and conjugacy class (labelled
by the order of elements in the class), where we write $x$ for $p^{-s}$.
For $S_5$ it turns out that the Frobenius conjugacy class at $p$ is
determined by the number of linear and quadratic factors of the reduction
$\bar{f}$ of $f$ modulo $p$, which may be computed from the degrees of
$\gcd\bigl(x^{p^n}-x,\bar{f}(x)\bigr)$ for $n=1,2$.  That computation
requires $O(\log p)$ multiplications and additions mod $p$.  Thus,
by the prime number theorem, for each $L(s,\rho)$ we may determine the
Dirichlet coefficients $a_n$ for $n\le X$ in time $O(X)$ (assuming mod
$p$ multiplications and additions take bounded time, which is appropriate
for numbers of the size that we consider).  Up to the implied constant,
that is best possible.  Moreover, the technique is very fast in practice;
we found that it takes approximately seven hours to compute the local factors
for all $p<2^{32}$.

\begin{table}[h]
\begin{tabular}{|lllll|}
\hline
$\rho$ & 1 & 2a & 2b & 3 \\ \hline
$1$ & $1-x$ & $1-x$ & $1-x$ & $1-x$ \\
$\chi$ & $1-x$ & $1+x$ & $1-x$ & $1-x$ \\
$\rho_4$ & $(1-x)^4$ & $(1-x)^2(1-x^2)$ &
$(1-x^2)^2$ & $(1-x)(1-x^3)$ \\
$\rho_4'$ & $(1-x)^4$ & $(1+x)^2(1-x^2)$ &
$(1-x^2)^2$ & $(1-x)(1-x^3)$ \\
$\rho_5$ & $(1-x)^5$ & $(1+x)(1-x^2)^2$ &
$(1-x)(1-x^2)^2$ & $(1+x+x^2)(1-x^3)$ \\
$\rho_5'$ & $(1-x)^5$ & $(1-x)(1-x^2)^2$ &
$(1-x)(1-x^2)^2$ & $(1+x+x^2)(1-x^3)$ \\
$\rho_6$ & $(1-x)^6$ & $(1-x^2)^3$ &
$(1+x)^2(1-x^2)^2$ & $(1-x^3)^2$ \\ \hline\hline
$\rho$ & 4 & 5 & 6 & \\ \hline
$1$ & $1-x$ & $1-x$ & $1-x$ & \\
$\chi$ & $1+x$ & $1-x$ & $1+x$ & \\
$\rho_4$ & $1-x^4$ &
$1+x+x^2+x^3+x^4$ & $(1+x)(1-x^3)$ & \\
$\rho_4'$ & $1-x^4$ &
$1+x+x^2+x^3+x^4$ & $(1-x)(1+x^3)$ & \\
$\rho_5$ & $(1-x)(1-x^4)$ &
$1-x^5$ & $(1+x+x^2)(1+x^3)$ & \\
$\rho_5'$ & $(1+x)(1-x^4)$ &
$1-x^5$ & $(1-x+x^2)(1-x^3)$ & \\
$\rho_6$ & $(1+x^2)(1-x^4)$ &
$(1-x)(1-x^5)$ & $1-x^6$ & \\ \hline
\end{tabular}
\caption{$S_5$ unramified local factors}
\label{tab:S5factors}
\end{table}

\subsection{Estimating zeros} \label{numerics:explicit}
\begin{figure}[h]
\centerline{
\epsfxsize=3in \epsfbox{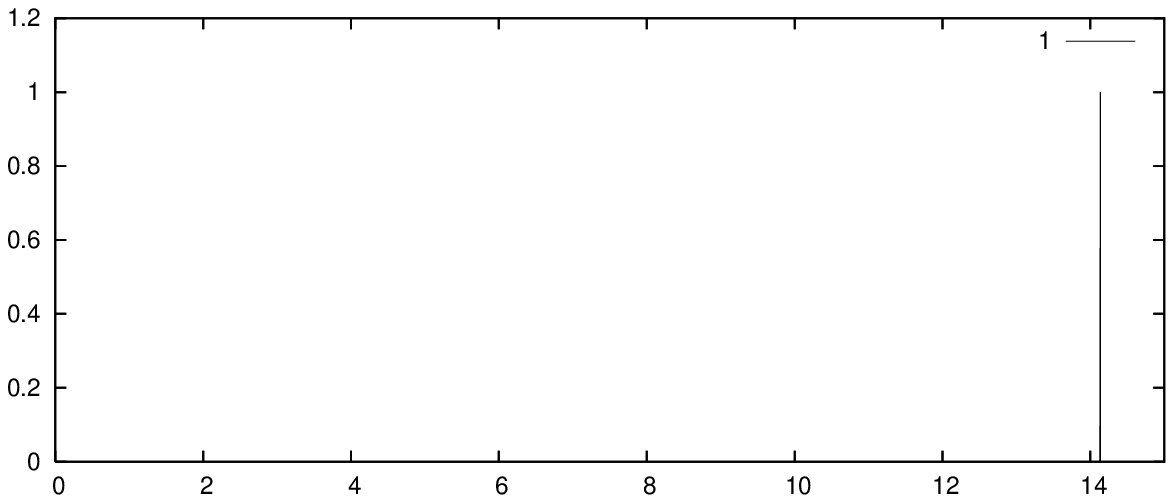}
\epsfxsize=3in \epsfbox{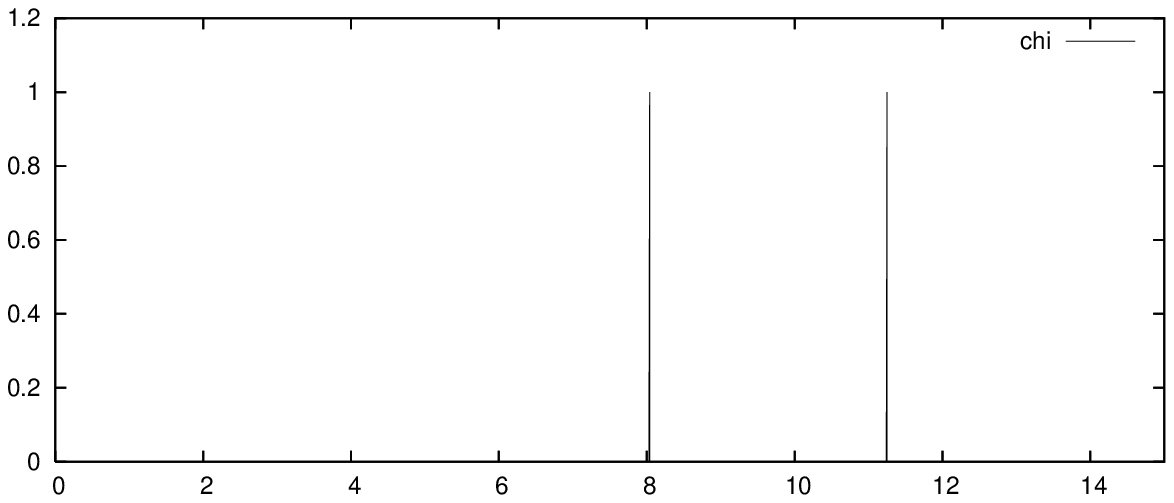}
}
\centerline{
\epsfxsize=3in \epsfbox{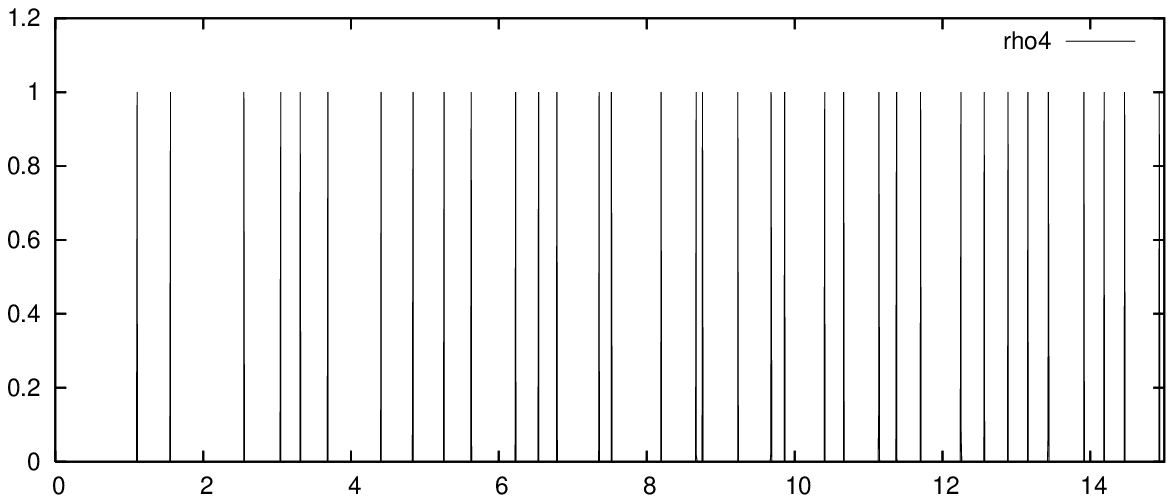}
\epsfxsize=3in \epsfbox{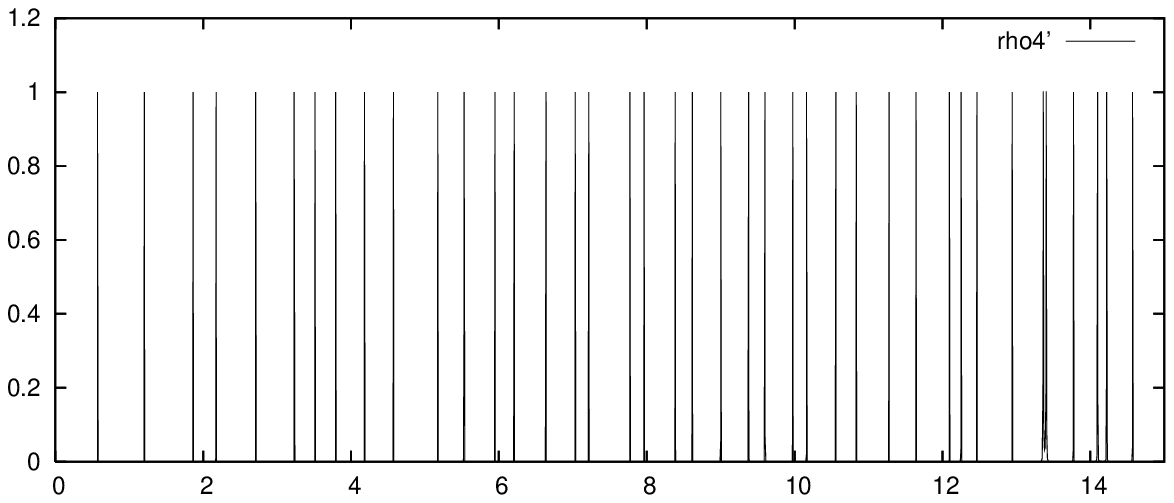}
}
\centerline{
\epsfxsize=3in \epsfbox{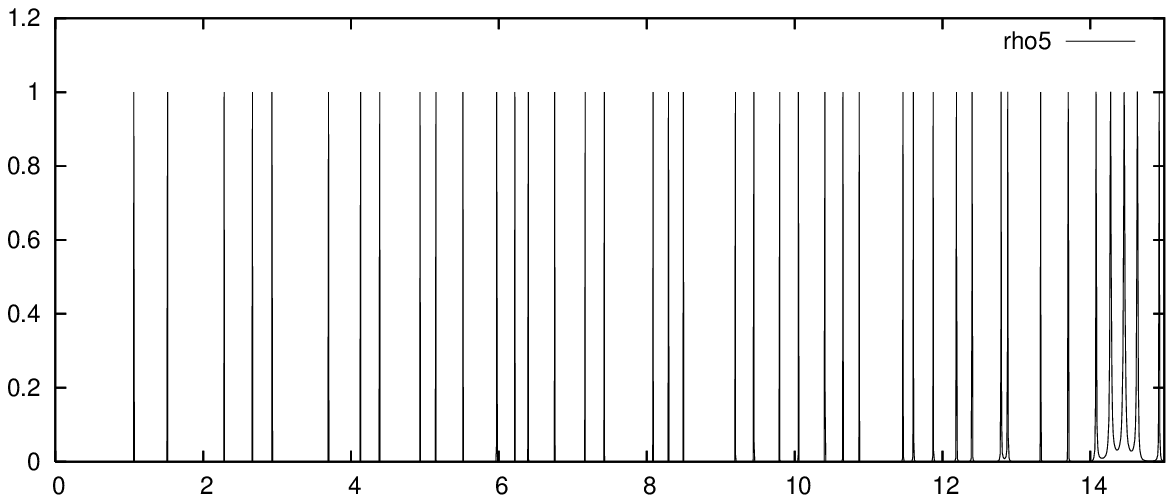}
\epsfxsize=3in \epsfbox{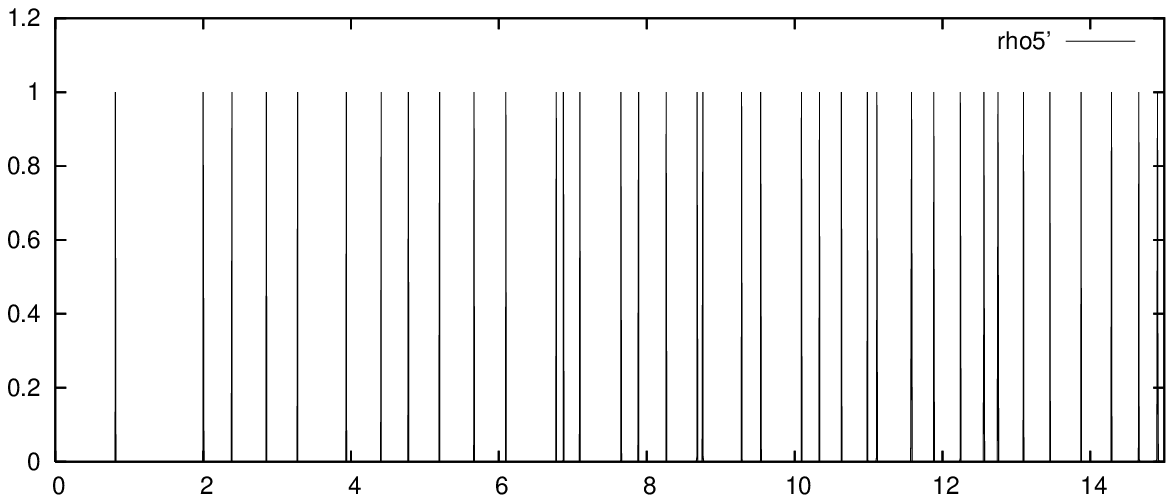}
}
\centerline{
\epsfxsize=6.4in \epsfbox{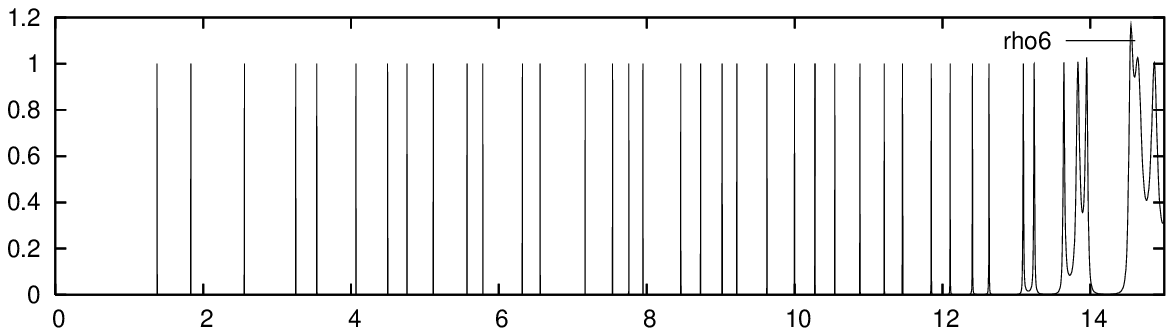}
}
\caption{$F_X(t)$ for each irreducible representation}
\label{fig:FXt}
\end{figure}

With our computed coefficients, we readily obtain estimates for the low
zeros by the method of Section \ref{sec:explicit}.  Figure \ref{fig:FXt}
shows graphs of $F_X(t)$, with $X=\log(2^{32})$, for each irreducible
$L$-function.  The spikes correspond to zeros, from which we get the
estimates for the ordinate of the lowest zero of each function shown in
Table \ref{tab:zeros}; note that for $\zeta$ the estimate agrees with
the known value $14.1347251417\ldots$ to within the precision of
the computation.  The increase in density of zeros with the conductor
and degree is apparent in the graphs.  Moreover, as the explicit formula
is very sensitive to errors in the coefficients, the fact that we see
spikes of height $1$ for the low zeros indicates that our coefficients
were computed correctly.  Each graph took a few minutes to generate.

\begin{table}
\begin{tabular}{|r|r|l||r|r|l|}
\hline
$\rho$ & conductor & lowest zero & $\rho$ & conductor & lowest zero \\ \hline
$1$       &        $1$ & $14.134725142$ &
$\chi$    &        $3$ & $8.039737156$\\
$\rho_4$  &  $4009008$ & $1.108937765$ &
$\rho_4'$ & $36081072$ & $0.5717508665$\\
$\rho_5$  & $36081072$ & $1.062064850$ &
$\rho_5'$ & $12027024$ & $0.8132800720$\\
$\rho_6$  & $36081072$ & $1.376872200$ &&& \\ \hline
\end{tabular}
\caption{Conductor and ordinate of the lowest zero of each
irreducible $L$-function}
\label{tab:zeros}
\end{table}

\subsection{Computing $G^{(k)}(u_m;\eta,\{\mu_j\})$}
\begin{figure}[h]
\centerline{\epsfxsize=6.1in \epsfbox{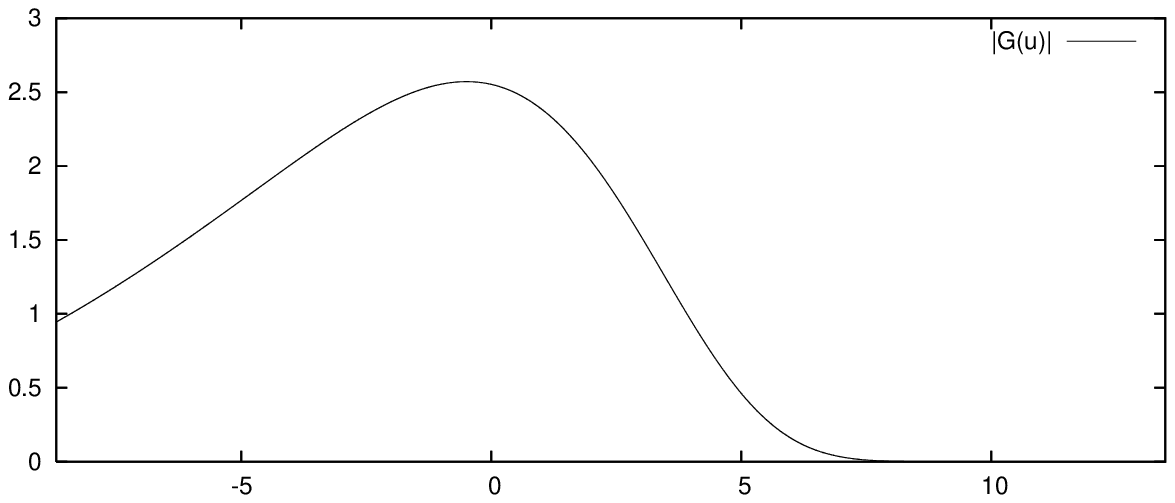}}
\caption{$|G(u;0.98,\{0,0,0,1,1,1\})|$}
\label{fig:Ggraph}
\end{figure}

Next we compute local approximations of $G(u;\eta,\{\mu_j\})$ for $u$
in the interval $\Bigl[\log\frac1{\sqrt N},\log\frac{2^{32}}{\sqrt N}\Bigr]$.
We evaluate $2^{13}$ Taylor series of $16$ terms using \eqref{powerlog};
with these choices, the error term in \eqref{localapprox} is less than
$10^{-28}$.  This calculation is the most delicate, due to high precision
and catastrophic cancellation.  Nevertheless, the computation time for
this stage was only a few hours.  The graph of $|G(u;\eta,\{\mu_j\})|$
is shown in Figure \ref{fig:Ggraph}.

\subsection{Computing $S^{(k)}_m$, \eqref{eq:Fxlocal} and $L(s,\rho_6)$}
Now we come to the main part of the computation, \eqref{eq:Fxlocal}.
Most of the time, approximately twelve hours, was spent computing
$S^{(k)}_m$.  Note that if we had not adjusted $\eta$ to reduce the
precision, this calculation could have taken substantially longer.

Once we have $S^{(k)}_m$, the computation of
\eqref{eq:Fxlocal} and $L(s,\rho_6)$ is very fast.  We choose
$B=\frac{2\pi\cdot2^{12}}{\log(2^{32})}\approx1160$.  Since this is much
larger than $t=115$, the errors terms from Lemma \ref{lem:Lambdasum} are
negligible.  We choose $A=2^{20}/B\approx900$, which is about $160$ times
the expected density $\frac1{2\pi}\log N\!\left(\frac{t}{2\pi}\right)^r$
of zeros around $t=115$.  Thus, the main Fourier transform is of $2^{20}$
points, which takes only a few minutes to compute.

\begin{figure}[h]
\centerline{\epsfxsize=6.1in \epsfbox{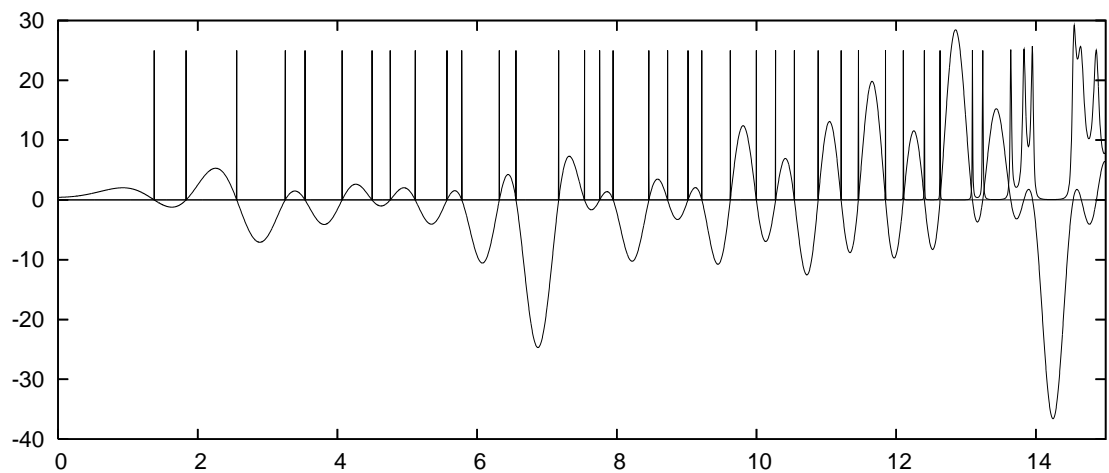}}
\caption{$Z(t)$ and $25F_X(t)$ for small $t$}
\label{fig:Lgraph1}
\end{figure}

\begin{figure}[h]
\centerline{\epsfxsize=6.1in \epsfbox{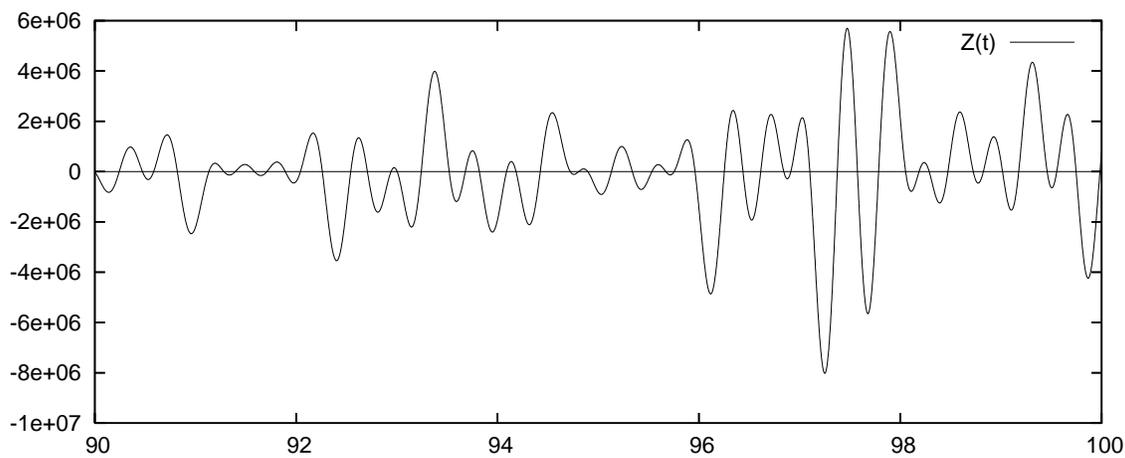}}
\caption{$Z(t)$ for large $t$}
\label{fig:Lgraph2}
\end{figure}

Figure \ref{fig:Lgraph1} shows the graph of
$Z(t):=\Lambda\bigl(\frac12+it\bigr)/\bigl|\gamma\bigl(\frac12+it\bigr)\bigr|$,
which is the analogue of Riemann's $Z$ function.  We have
superimposed the graph of $25F_X(t)$ over the same range; note the good
agreement in location of zeros between the two, which gives evidence
that our computations are correct.  Figure \ref{fig:Lgraph2} shows
$Z(t)$ over the higher range $t\in [90,100]$.

\subsection{Turing's method}
Finally, we apply Turing's method to the computed $L$-functions.  At the
same time, we verify the ``working hypothesis'' that we can isolate the
zeros of the irreducible $L$-functions.  This verification takes only a
few seconds.

\clearpage
\bibliographystyle{alpha}
\bibliography{turing}

\begin{thebibliography}{BDSBT01}

\bibitem[Arm72]{armitage}
J.V. Armitage.
\newblock Zeta functions with a zero at $s=\frac12$.
\newblock {\em Invent. Math.}, 16:195--205, 1972.

\bibitem[Art30]{artin}
E.~Artin.
\newblock Zur {T}heorie der {L}-{R}eihen mit allgemeinen {G}ruppencharakteren.
\newblock {\em Abh. Math. Sem. Univ. Hamburg}, 8:292--306, 1930.

\bibitem[BDSBT01]{taylor1}
K.~Buzzard, M.~Dickinson, N.~Shepard-Barron, and R.~Taylor.
\newblock On icosahedral {A}rtin representations.
\newblock {\em Duke Math. J.}, 109:283--318, 2001.

\bibitem[Boo03a]{booker2}
Andrew~R. Booker.
\newblock {\em Numerical tests of modularity}.
\newblock PhD thesis, Princeton University, 2003.
\newblock available from {\tt http://www.umich.edu/$\sim$arbooker/papers/}.

\bibitem[Boo03b]{booker}
Andrew~R. Booker.
\newblock Poles of {A}rtin {$L$}-functions and the strong {A}rtin conjecture.
\newblock {\em Ann. of Math. (2)}, 158(3):1089--1098, 2003.

\bibitem[Bra47]{brauer}
R.~Brauer.
\newblock On {A}rtin's {$L$}-series with general group characters.
\newblock {\em Ann. Math.}, 48:502--514, 1947.

\bibitem[BS02]{stein}
K.~Buzzard and W.~Stein.
\newblock A mod five approach to modularity of icosahedral {G}alois
  representations.
\newblock {\em Pac. Jour. Math.}, 203(2):265--282, 2002.

\bibitem[BS05]{bs}
A.~Booker and A.~Str\"ombergsson.
\newblock Numerical computations with the trace formula and the {S}elberg
  eigenvalue conjecture.
\newblock {\em in preparation}, 2005.

\bibitem[BSV05]{bsv}
A.~Booker, A.~Str\"ombergsson, and A.~Venkatesh.
\newblock Effective computation of {M}aass cusp forms.
\newblock {\em preprint}, 2005.

\bibitem[Buh78]{buhler}
J.P. Buhler.
\newblock {\em Icosahedral {G}alois representations}.
\newblock Springer-Verlag, Berlin, 1978.
\newblock Lecture Notes in Mathematics, Vol. 654.

\bibitem[Dok04]{dokchitser}
Tim Dokchitser.
\newblock Computing special values of motivic {$L$}-functions.
\newblock {\em Experiment. Math.}, 13(2):137--149, 2004.

\bibitem[DS74]{deligne-serre}
Pierre Deligne and Jean-Pierre Serre.
\newblock Formes modulaires de poids {$1$}.
\newblock {\em Ann. Sci. \'Ecole Norm. Sup. (4)}, 7:507--530 (1975), 1974.

\bibitem[Dum91]{dummit}
D.~S. Dummit.
\newblock Solving solvable quintics.
\newblock {\em Math. Comp.}, 57(195):387--401, 1991.

\bibitem[Fli94]{flicker}
Yuval~Z. Flicker.
\newblock On the symmetric square: total global comparison.
\newblock {\em J. Funct. Anal.}, 122(2):255--278, 1994.

\bibitem[FM89]{foote-murty}
Richard Foote and V.~Kumar Murty.
\newblock Zeros and poles of {A}rtin {$L$}-series.
\newblock {\em Math. Proc. Cambridge Philos. Soc.}, 105(1):5--11, 1989.

\bibitem[GAP05]{GAP}
The GAP~Group.
\newblock {\em {GAP -- Groups, Algorithms, and Programming, Version 4.4}},
  2005.
\newblock \verb+(http://www.gap-system.org)+.

\bibitem[IS00]{iwaniec-sarnak}
H.~Iwaniec and P.~Sarnak.
\newblock Perspectives on the analytic theory of {$L$}-functions.
\newblock {\em Geom. Funct. Anal.}, (Special Volume, Part II):705--741, 2000.
\newblock GAFA 2000 (Tel Aviv, 1999).

\bibitem[JM00]{jehanne}
A.~Jehanne and M.~M{\"u}ller.
\newblock Modularity of an odd icosahedral representation.
\newblock {\em J. Th\'eor. Nombres Bordeaux}, 12(2):475--482, 2000.
\newblock Colloque International de Th\'eorie des Nombres (Talence, 1999).

\bibitem[JM01]{jehanne2}
A.~Jehanne and M.~M{\"u}ller.
\newblock Modularity of some odd icosahedral representations.
\newblock {\em unpublished}, 2001.
\newblock available from {\tt
  http://www.math.u-bordeaux1.fr/$\sim$jehanne/travaux.html}.

\bibitem[Kim94]{kiming}
I.~Kiming.
\newblock On the experimental verification of the {A}rtin conjecture for
  {$2$}-dimensional odd {G}alois representations over {${\bf Q}$}. {L}iftings
  of {$2$}-dimensional projective {G}alois representations over {${\bf Q}$}.
\newblock In {\em On Artin's conjecture for odd $2$-dimensional
  representations}, volume 1585 of {\em Lecture Notes in Math.}, pages 1--36.
  Springer, Berlin, 1994.

\bibitem[KM]{nfdatabase}
J\"urgen Kl\"uners and Gunter Malle.
\newblock A database for number fields.
\newblock available from {\tt
  http://www.mathematik.uni-kassel.de/$\sim$klueners/minimum/minimum.html}.

\bibitem[Lan80]{langlands}
R.P. Langlands.
\newblock {\em Base change for ${\rm {G}{L}}(2)$}.
\newblock Princeton University Press, Princeton, N.J., 1980.

\bibitem[Leh70]{lehman}
R.~Sherman Lehman.
\newblock On the distribution of zeros of the {R}iemann zeta-function.
\newblock {\em Proc. London Math. Soc. (3)}, 20:303--320, 1970.

\bibitem[LO79]{lagarias-odlyzko}
J.~C. Lagarias and A.~M. Odlyzko.
\newblock On computing {A}rtin {$L$}-functions in the critical strip.
\newblock {\em Math. Comp.}, 33(147):1081--1095, 1979.

\bibitem[LRS99]{lrs}
Wenzhi Luo, Ze{\'e}v Rudnick, and Peter Sarnak.
\newblock On the generalized {R}amanujan conjecture for {${\rm GL}(n)$}.
\newblock In {\em Automorphic forms, automorphic representations, and
  arithmetic (Fort Worth, TX, 1996)}, volume~66 of {\em Proc. Sympos. Pure
  Math.}, pages 301--310. Amer. Math. Soc., Providence, RI, 1999.

\bibitem[Odl87]{odlyzko}
A.~M. Odlyzko.
\newblock On the distribution of spacings between zeros of the zeta function.
\newblock {\em Math. Comp.}, 48:273--308, 1987.

\bibitem[Oma01]{omar}
Sami Omar.
\newblock Localization of the first zero of the {D}edekind zeta function.
\newblock {\em Math. Comp.}, 70(236):1607--1616 (electronic), 2001.

\bibitem[OS88]{odlyzko-schonhage}
A.~M. Odlyzko and A.~Sch{\"o}nhage.
\newblock Fast algorithms for multiple evaluations of the {R}iemann zeta
  function.
\newblock {\em Trans. Amer. Math. Soc.}, 309(2):797--809, 1988.

\bibitem[RR05]{MPFI}
Nathalie Revol and Fabrice Rouiller.
\newblock Multiple precision floating-point interval library, version {\tt
  1.3.3}, 2005.
\newblock available from {\tt
  http://perso.ens-lyon.fr/nathalie.revol/software.html}.

\bibitem[RS96]{rudnick-sarnak}
Ze{\'e}v Rudnick and Peter Sarnak.
\newblock Zeros of principal {$L$}-functions and random matrix theory.
\newblock {\em Duke Math. J.}, 81(2):269--322, 1996.
\newblock A celebration of John F. Nash, Jr.

\bibitem[Rub05]{rubinstein}
Michael Rubinstein.
\newblock Computational methods and experiments in analytic number theory.
\newblock {\em preprint}, 2005.

\bibitem[Rum93]{rumely}
Robert Rumely.
\newblock Numerical computations concerning the {ERH}.
\newblock {\em Math. Comp.}, 61(203):415--440, S17--S23, 1993.

\bibitem[Tay03]{taylor2}
Richard Taylor.
\newblock On icosahedral {A}rtin representations. {II}.
\newblock {\em Amer. J. Math.}, 125(3):549--566, 2003.

\bibitem[{The}04]{PARI}
{The PARI~Group}, Bordeaux.
\newblock {\em {PARI/GP, version {\tt 2.1.5}}}, 2004.
\newblock available from {\tt http://pari.math.u-bordeaux.fr/}.

\bibitem[Tit86]{titchmarsh}
E.~C. Titchmarsh.
\newblock {\em The theory of the {R}iemann zeta-function}.
\newblock The Clarendon Press Oxford University Press, New York, second
  edition, 1986.
\newblock Edited and with a preface by D. R. Heath-Brown.

\bibitem[Tol97]{tollis}
Emmanuel Tollis.
\newblock Zeros of {D}edekind zeta functions in the critical strip.
\newblock {\em Math. Comp.}, 66(219):1295--1321, 1997.

\bibitem[Tun81]{tunnell}
J.~Tunnell.
\newblock Artin's conjecture for representations of octahedral type.
\newblock {\em Bull. AMS}, 5:173--175, 1981.

\bibitem[Tur53]{turing}
A.~M. Turing.
\newblock Some calculations of the {R}iemann zeta-function.
\newblock {\em Proc. London Math. Soc. (3)}, 3:99--117, 1953.

\bibitem[Wan03]{wang}
Song Wang.
\newblock On the symmetric powers of cusp forms on {${\rm GL}(2)$} of
  icosahedral type.
\newblock {\em Int. Math. Res. Not.}, (44):2373--2390, 2003.

\bibitem[Wei99]{weissman}
M.~Weissman.
\newblock Icosahedral {G}alois representations and modular forms.
\newblock {\em Princeton University}, Undergraduate thesis, 1999.
\newblock available from {\tt http://math.berkeley.edu/$\sim$marty/}.

\end{thebibliography}
\end{document}